\documentclass[12pt,a4paper]{article}
\usepackage{amsfonts}
\usepackage{amsmath}
\usepackage{graphicx}
\usepackage{txfonts}

\newtheorem{theorem}{Theorem}

\newtheorem{pr}[theorem]{Proposition}
\newtheorem{re}[theorem]{Remark}

\DeclareMathOperator{\sech}{sech}
 \DeclareMathOperator{\sn}{sn}
\DeclareMathOperator{\dn}{dn} \DeclareMathOperator{\cn}{cn}
\date{}
\begin{document}

\title{\textbf{Mathematical pendulum and its variants}}
\author{ \textbf{O. Chi\c{s}$^*$, D. Opri\c{s}$^{**}$}}
\maketitle
\begin{center}
$^*, \, ^{**}$  Faculty of Mathematics and Informatics, West University of Timi\c{s}oara, Romania\\
E-mail: chisoana@yahoo.com, opris@math.uvt.ro
\end{center}

\textbf{Abstract}: In this paper we show that there are
applications that transform the movement of a pendulum into
movements in $\mathbb{R}^3$. This can be done using Euler top
system of differential equations. On the constant level surfaces,
Euler top system reduces to the equation of a pendulum. Those
properties are also considered in the case of system of
differential equations with delay argument and in the fractional
case. A\-nother aspect presented here is stochastic Euler top
system of
differential equations and stochastic pendulum.\\

\maketitle \textbf{MSC2000}: 34K50, 35L65, 26A33, 37N99, 60H10, 65C20, 65C30.\\


\section{Introduction}

\par The dynamics of some mechanical systems is described using the
rigid body dynamics with a fixed point, mathematical pendulum or
oscillators. These systems belong to a class of differential
equations from  $\mathbb{R}^3$ with the right side polynomial
functions of degree greater or equal to two. From this category we
will consider Euler top system of differential equations. We begin
our study from mathematical pendulum (and its variants: with
delay, fractional and stochastic) approach.

The Euler top system of fractional differential equation belongs
to a class of differential equations that are described using
polynomial functions. It has the form
 \begin{equation}\label{1}
\left\{%
\begin{array}{ll}
    \dot{x_{1}}(t)=x_{2}(t)x_{3}(t), &  \\
    \dot{x_{2}}(t)=-x_{1}(t)x_{3}(t), &  \\
    \dot{x_{3}}(t)=x_{1}(t)x_{2}(t).&  \\
\end{array}%
\right.
    \end{equation}

Because system (\ref{1}) has three Hamilton-Poisson realizations,
three conservation laws are given by the Hamiltonians $H_1, \,
H_2$ and $H_3$ \cite{Ch2}:
\begin{enumerate}
    \item $H_{1}(x_{1}(t), x_{2}(t), x_{3}(t)):= \frac{1}{2}(x_{1}^{2}(t)+x_{2}^{2}(t));$
    \item $H_{2}(x_{1}(t), x_{2}(t), x_{3}(t)):= -\frac{1}{2}(x_{2}^{2}(t)+x_{3}^{2}(t));$
    \item $H_{3}(x_{1}(t), x_{2}(t), x_{3}(t)):= x_{1}^{2}(t)-x_{3}^{2}(t),$
\end{enumerate}
and the other three conservation laws are given by the
corresponding Casimir functions of the above realizations
\cite{Ch2}:
\begin{enumerate}
    \item $C_{1}(x_{1}(t), x_{2}(t),
x_{3}(t)):=\frac{1}{2}(x_{2}^{2}(t)+x_{3}^{2}(t));$
    \item $C_{2}(x_{1}(t), x_{2}(t),
x_{3}(t)):=\frac{1}{2}(x_{1}^{2}(t)-x_{2}^{2}(t));$
    \item $C_{3}(x_{1}(t), x_{2}(t),
x_{3}(t)):=x_{1}^{2}(t)+x_{2}^{2}(t).$
\end{enumerate}

A simple mathematical pendulum is the mathematical model of a
ball, having the mass $m,$ which hangs in a point  $O$ by a bar of
length  $l,$ and the point  $O$ performs movement in a plane
\cite{Trueba}.

The Euler-Lagrange equation that describes the movement of a
pendulum is given by
\begin{equation}\label{Pendul}
l\ddot{\theta}(t)+g\sin\theta(t)+\ddot{x}_0\cos\theta(t)-\ddot{y}_0\sin\theta(t)=0.
\end{equation}
The dumping pendulum equations with periodic force is
\begin{equation}\label{pendul2}
\ddot{\theta}(t)+2h\sin\theta(t)+f_1(t)\cos\theta(t)+f_2(t)\sin\theta(t)+
\sum_{p=0}^N \alpha_p\dot{\theta}(t)|\dot{\theta}(t)|^{p-1}=0,
\end{equation}
and for $f_1:=0, \, f_2:=0$ and $\alpha_p:=0, \, p=1...N,$ then
(\ref{pendul2}) reduces to $\ddot{\theta}(t)+2h\sin\theta(t)=0.$

In the first section we will determine the analytical solutions
for Euler top system taking into consideration the conservation
laws that it owns, and point out the analytical solution for
pendulum. In the second section we have presented the co\-nnection
between Euler top system and pendulum: the restriction of the
system to a constant level surface represents the pendulum
equations. The third section presents the Euler top system of
differential equations with delay argument, along the OZ and OX
axes. These new systems have also conservation laws and the
restriction of the orbits at these surfaces of constant level
determined by the conservation laws are mathematical pendulums
with delay argument. In the forth section we presented the Euler
top system of fractional differential equations. We have used
Caputo fractional derivative in OZ and OX directions. As in the
previous case, this system of fractional differential equations
have conservation laws and the restriction of the system to the
constant level surfaces is a fractional pendulum. In Section 5 we
presented stochastic Euler top system and stochastic pendulum. We
considered It\^o and Stratonovich integrals for describing the
stochastic process, using a Wiener process. For all these cases
numerical simulations are done. In the last section some
conclusions are presented and ideas for future work.

\section{Euler top system and simple pendulum - analytical solutions}

Let us consider the Euler top system of differential equations
(\ref{1}) and the integrals of motion given by
\begin{equation}\label{2}
x_{1}^{2}(t)+x_{2}^{2}(t)=2H^2, \quad
x_{2}^{2}(t)+x_{3}^{2}(t)=2K^2.
\end{equation}

From (\ref{2}), results that
\begin{equation}\label{3}
x_{1}^{2}(t)=2H^2-x_{2}^{2}(t), \quad
x_{3}^{2}(t)=2K^2-x_{2}^{2}(t).
\end{equation}
Replacing (\ref{3}) in the first equation in (\ref{1}) we get:
\begin{equation}
(\dot{x}_2)^2(t)=(x_1)^2+(x_3)^2(t)=(2H^2-x_{2}^{2}(t))(2K^2-x_{3}^{2}(t))
\end{equation}
and so,
\begin{equation}\label{5}
t=\int_{x_2(0)}^{x_2(t)}\frac{1}{\sqrt{(2H^2-u^2)(2K^2-u^2)}}du,
\end{equation}
that shows that $x_2(t),$ respectively $x_1(t)$ and $x_2(t)$ are
elliptic functions of time \cite{MR}.

In the case when $H=K,$ the quartic under the square root has
double roots and (\ref{5}) can be explicitly integrated by means
of elementary functions in the following manner. The equation
$$\dot{x}_2(t)=\pm (2H^2-x_2^2(t)),$$ with $x_2(0)=0,$ has the solution
\begin{equation}\label{sol2}
x_2(t)=\pm H\sqrt{2}\tanh(H\sqrt{2}t).
\end{equation}
Substituting (\ref{sol2}) in (\ref{5}), we get
\begin{equation}\label{sol1+3}
x_1(t)=\pm H\sqrt{2}\sech(H\sqrt{2}t), \quad x_3(t)=\pm
H\sqrt{2}\sech(H\sqrt{2}t).
\end{equation}
So, the equations (\ref{sol1+3}) and (\ref{pendul2}) represent the
two heteroclinic orbits for the Euler top system and are given by
$$\Big(\pm H\sqrt{2}\sech(H\sqrt{2}t), \pm H\sqrt{2}\tanh(H\sqrt{2}t),
\pm H\sqrt{2}\sech(H\sqrt{2}t)\Big).$$

In the case when $H\neq K,$ the integral (\ref{5}) can be computed
using Jacobi's elliptic functions \cite{Law}. We use relations
$$\frac{d}{dt}\sn u=\cn u \dn u,\quad \cn^2 u=1-sn^2u, \quad \dn^2 u=1-m^2
\sn^2 u$$ and
\begin{equation}\label{sola2}
x_2(t)=H\sqrt{2}\sn
\Big(H\sqrt{2}t;\frac{\sqrt{H}}{\sqrt{K}}\Big),
\end{equation}
 with the initial condition
$x_2(0)=0.$ Choosing the time deviation, appropriately, we can
assume that $\dot{x}_2(0)>0.$ From (\ref{3}) results that
\begin{equation}\label{sola1+3}
x_1(t)=H\sqrt{2}\cn
\Big(H\sqrt{2}t;\frac{\sqrt{H}}{\sqrt{K}}\Big),\quad
x_3(t)=K\sqrt{2}\sn
\Big(H\sqrt{2}t;\frac{\sqrt{H}}{\sqrt{K}}\Big).
\end{equation}
If $\phi$ denotes the period invariant of Jacobi's elliptic
functions, then $x_1(t)$ and $x_2(t)$ have the period
$4\phi/H\sqrt{2},$ whereas $x_3(t)$ has the period
$2\phi/H\sqrt{2}.$
\begin{pr}
\begin{description}
    \item[a)] If $H=K,$ than Euler top system \emph{(\ref{1})} has an
    analytical solution given by \emph{(\ref{sol2})} and \emph{(\ref{sol1+3})};
    \item[b)] If $H\neq K,$ than the Euler top system has the
    analytical solution given by \emph{(\ref{sola2})} \\ and
    \emph{(\ref{sola1+3})}.
\end{description}
\end{pr}\hfill $\Box$

\begin{pr}\emph{\cite{Bele}}
The analytical solution for simple pendulum
$\ddot{\theta}(t)+2h\sin\theta(t)=0,$ with initial conditions
$\theta(0)=\theta_0$ and $\dot{\theta}(0)=0$ is given by
$$\theta(t)=2\arcsin\Big\{\sin\theta_0 \sn\Big(\sin^2\frac{\theta_0}{2}\Big)-
\sqrt{2h}t;\sin^2\frac{\theta_0}{2}\Big\}.$$
\end{pr}\hfill $\Box$

\section{Euler top system and simple pendulum}

In this section we will show the way the Euler top system and the
simple pendulum are linked. We will show that the movement of the
Euler top system is reduced to pendulum movement on the constant
level surfaces $H$ and $K,$ described by the conservation laws:
\begin{equation}\label{c1}
\frac{1}{2}(x_1(t))^2+\frac{1}{2}(x_2(t))^2=H,
\end{equation}
\begin{equation}\label{c2}
\frac{1}{2}(x_2(t))^2+\frac{1}{2}(x_3(t))^2=K.
\end{equation}

Since $H$ and $K$ are conserved, the Euler top motion takes place
along the intersections of the level surfaces of the energy and
the angular momentum in $\mathbb{R}^3.$

\begin{pr}\label{pr1}
Let us consider the Euler top system of differential equations
\emph{(\ref{1})}.\begin{enumerate} \item The function $H,$ given
by \emph{(\ref{c1})} is a conservation law for system
\emph{(\ref{1})};
    \item The solution of \emph{(\ref{1})} on the constant level surface defined by
    \emph{(\ref{c1})}, given by
\begin{equation}
    (x_1(t))^2+(x_2(t))^2=2H=const, \, H>0
\end{equation} is
\begin{equation}\label{solutia1}
x_1(t)=\sqrt{2H}\cos \frac{\theta(t)}{2}, \quad
x_2(t)=\sqrt{2H}\sin \frac{\theta(t)}{2}, \quad
x_3(t)=-\frac{1}{2}\dot{\theta}(t),
\end{equation}where $\theta(t)$ is the solution of pendulum
equation $\ddot{\theta}(t)+2H\sin\theta(t)=0.$
\end{enumerate}
\end{pr}

\textbf{Proof:}\\
\emph{1.} By deriving (\ref{c1}) and by replacing it in (\ref{1}),
we have $$\dot{H}(t)=x_1(t)\dot{x}_1(t)+x_2(t)\dot{x}_2(t)=0.$$
And
so  $H$ is a conservation law.\\

 \emph{2.} Using a direct calculus, it can easily checked that (\ref{solutia1}) is a solution
 for (\ref{1}) and reciprocal.\hfill $\Box$

\begin{pr}\label{pr2}
\begin{enumerate}
\item The function $K,$ given by \emph{(\ref{c2})} is a
conservation law for the Euler top system \emph{(\ref{1})};
    \item The solution of \emph{(\ref{1})}, on the constant level surface
    defined by \emph{(\ref{c2})}, given by
\begin{equation}
    (x_2(t))^2+(x_3(t))^2=2K=const, \, K>0
\end{equation} is
\begin{equation}\label{solutia2}
x_1(t)=-\frac{1}{2}\dot{\theta}(t), \quad
x_2(t)=\sqrt{2K}\cos\frac{\theta(t)}{2}, \quad
x_3(t)=\sqrt{2K}\sin\frac{\theta(t)}{2},
\end{equation}where $\theta(t)$ is the solution of pendulum
equation $\ddot{\theta}(t)+2K\sin\theta(t)=0.$
\end{enumerate}
\end{pr}

\textbf{Proof:}\\
\emph{1.} By deriving (\ref{c2}) and by replacing it in (\ref{1}),
we have that $K$ is a conservation law because
$$\dot{K}(t)=x_2(t)\dot{x}_2(t)+x_3(t)\dot{x}_3(t)=0.$$
 \emph{2.} By direct calculations, it can be easily checked that
 (\ref{solutia2}) is a solution
 for (\ref{1}) and reciprocal.\hfill $\Box$

\begin{re}
The dynamics of Euler top system of differential equations in
$\mathbb{R}^3$ is a union of two-dimensional simple pendula.
\end{re}
\hfill $\Box$

For the initial conditions $x_1(0)=0.1, \, x_2(0)=0.1$ and
$x_3(0)=0.2,$ the Euler top system is represented in the first
figure and the pendulum is represented for the initial condition
$\theta(0)=-3.8.$
\begin{center}\begin{tabular}{ccc}
  \includegraphics[height=2.5in]{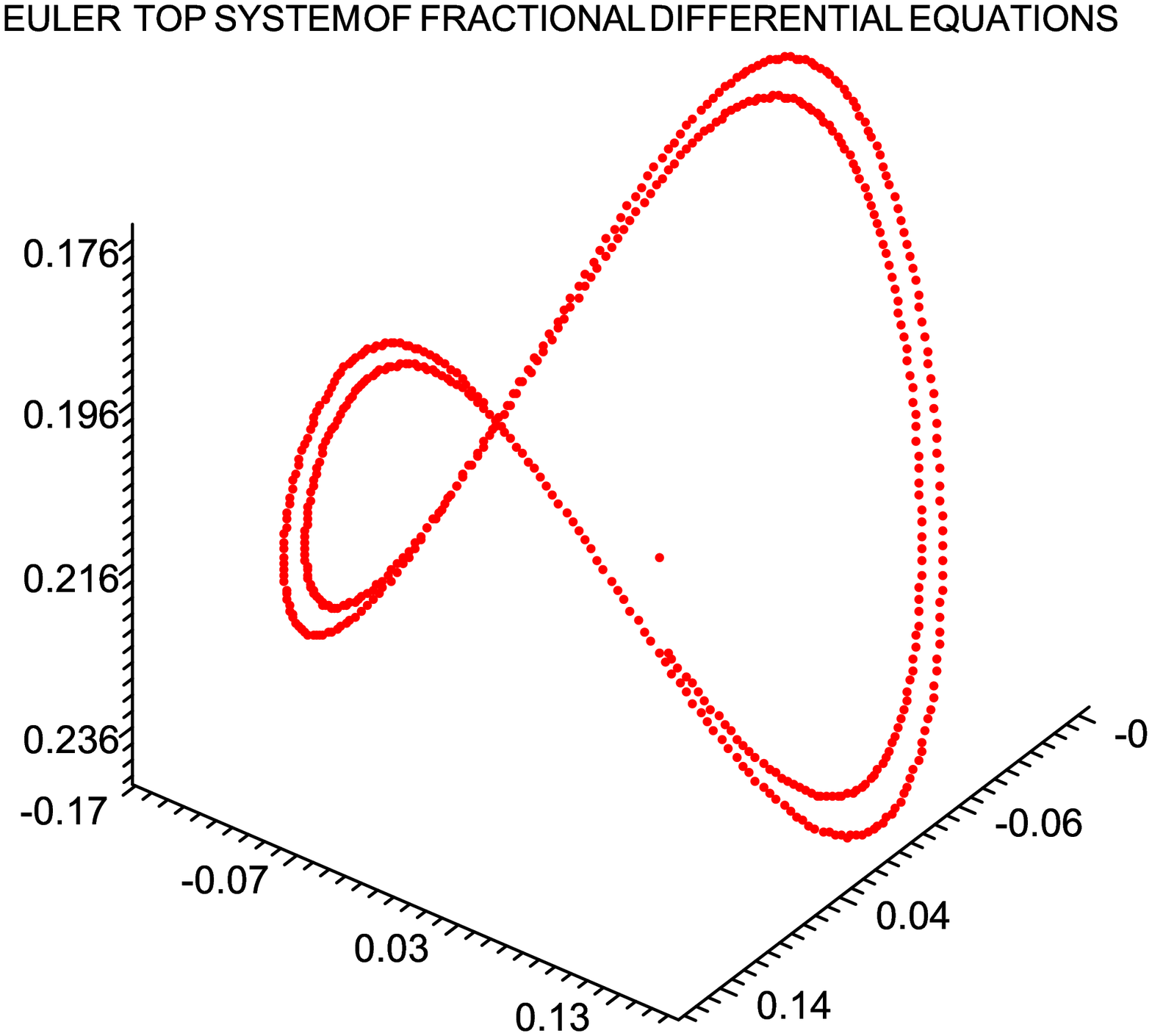}
  \includegraphics[height=2.5in]{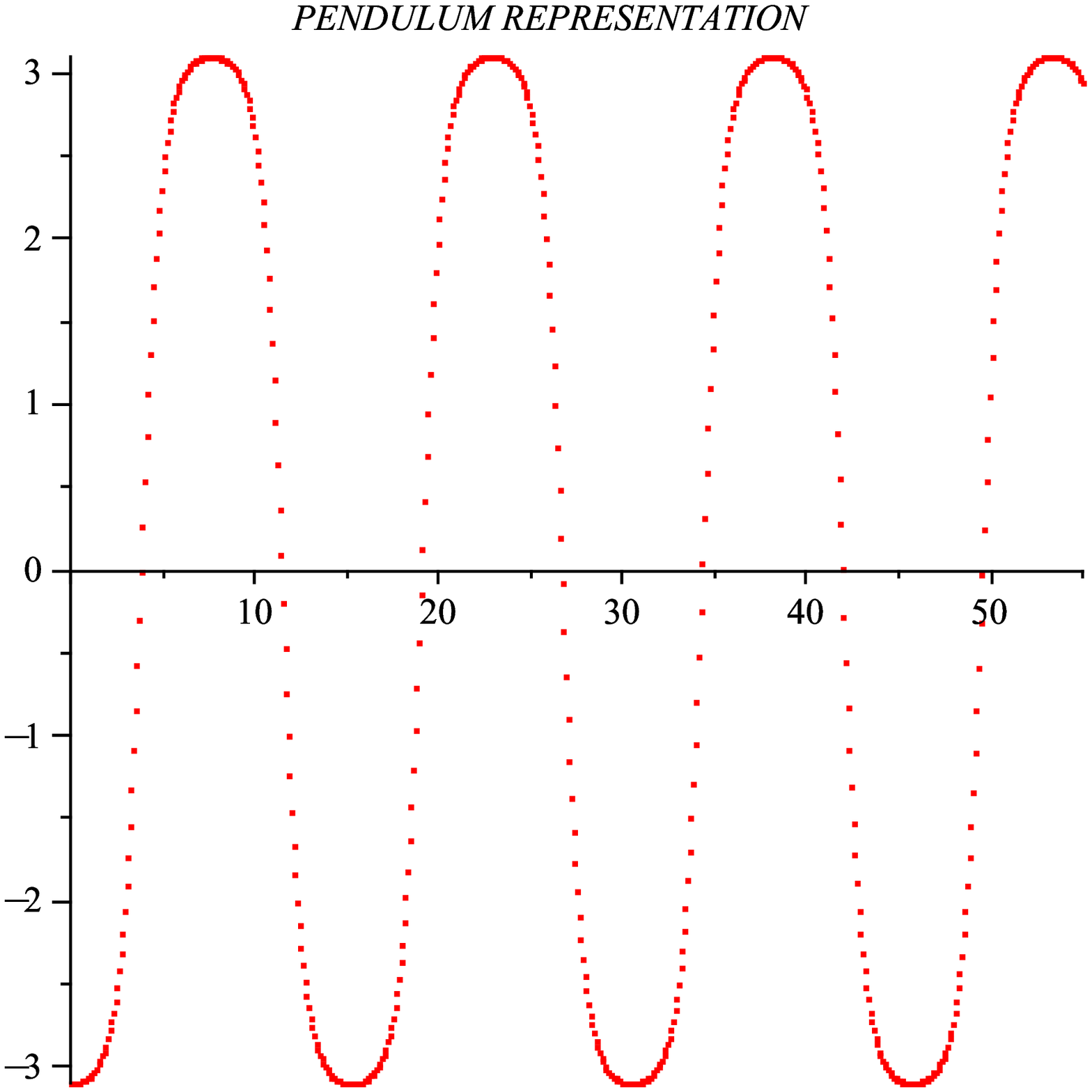}
\end{tabular}
\end{center}

\section{Euler top system and simple pendulum - with delay argument
and fractional derivative}

A differential equation with delay argument is defined in
\cite{Halle}. A second order diffe\-rential equation with delay
argument is given by
\begin{equation}\label{hale}
\ddot{\theta}(t)=c\sin(\theta(t-\tau)),
\end{equation}where $c\in\mathbb{R}$ is a solution of a
differential equation on the circle
$S^1=\{y\in\mathbb{R}^2|y_1^2+y_2^2=1\},$  $\theta$ an angle
variable determined up to a multiple of $2\pi,$ and $\tau>0.$

From Proposition \ref{pr1} and Proposition \ref{pr2} we can deduce
the following results.
\begin{pr}
The Euler top system of differential equations with delay argument
is given by
 \begin{equation}\label{1d}
\left\{%
\begin{array}{ll}
    \dot{x_{1}}(t)=x_{2}(t)x_{3}(t), &  \\
    \dot{x_{2}}(t)=-x_{1}(t)x_{3}(t), &  \\
    \dot{x_{3}}(t)=x_{1}(t-\tau)x_{2}(t-\tau).&  \\
\end{array}%
\right.
    \end{equation}
The system \emph{(\ref{1d})} has the following
properties\begin{description}
    \item[a)] The function $H$ given by \emph{(\ref{c1})};
    \item[b)] The solution of system \emph{(\ref{1d})} on the constant
    level surface \emph{(\ref{c1})} is given by \emph{(\ref{solutia1})}, where
    $\theta(t)$ is the solution of
    \begin{equation}
        \ddot{\theta}(t)+2H\sin\theta(t-\tau)=0
    \end{equation} and reciprocal.\hfill $\Box$
\end{description}
\end{pr}

\begin{pr}
The Euler top system of differential equations with delay argument
given by
 \begin{equation}\label{1dd}
\left\{%
\begin{array}{ll}
    \dot{x_{1}}(t)=x_{2}(t-\tau)x_{3}(t-\tau), &  \\
    \dot{x_{2}}(t)=-x_{1}(t)x_{3}(t), &  \\
    \dot{x_{3}}(t)=x_{1}(t)x_{2}(t),&  \\
\end{array}%
\right.
    \end{equation}
has the following properties
\begin{description}
    \item[a)] The function $K$ given by \emph{(\ref{c2})} is a
    conservation law for the system \emph{(\ref{1dd})};
    \item[b)] The solution of system \emph{(\ref{1dd})} on the
    constant level surface \emph{(\ref{c2})} is given by
    \emph{(\ref{solutia2})}, where $\theta(t)$ is solution of the
    equation
    \begin{equation}
        \ddot{\theta}(t)+2K\sin\theta(t-\tau)=0.
    \end{equation}
\end{description}
\end{pr}\hfill $\Box$

This system is considered to be a starting point in studying
differential equations with delay argument for differential
manifold.

For $H=0.5,$ and $\tau=1,$ the pendulum equation with delay
argument and with initial condition $\theta(0)=2,$ is represented
in the following figure. The Euler top system with delay argument
(\ref{1d}) is represented in the second figure, for the initial
conditions $x_1(0)=0.1, \, x_2(0)=0.05, \, x_3(0)=0.2.$
\begin{center}\begin{tabular}{cc}
  \includegraphics[height=2in]{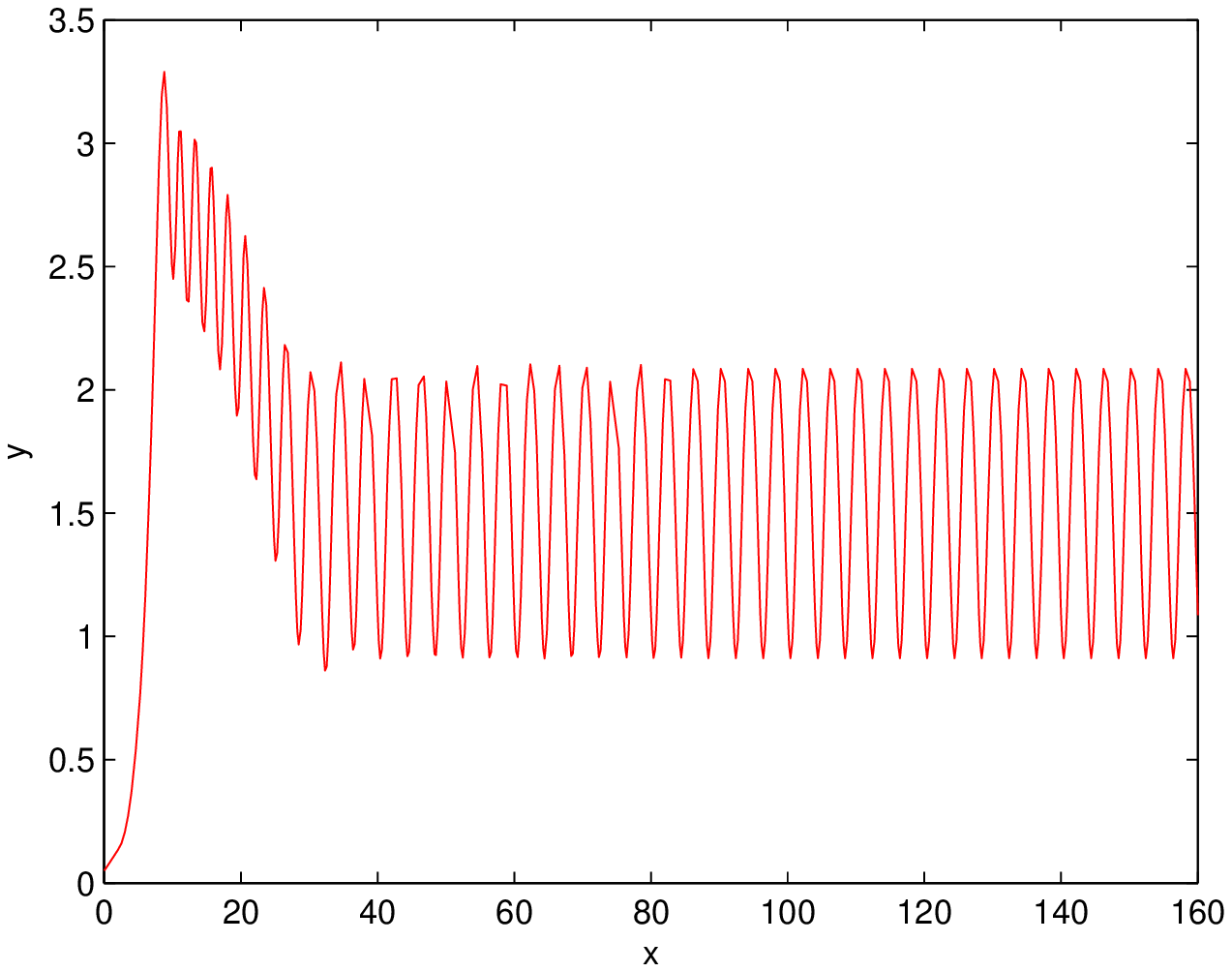} &
  \includegraphics[height=2in]{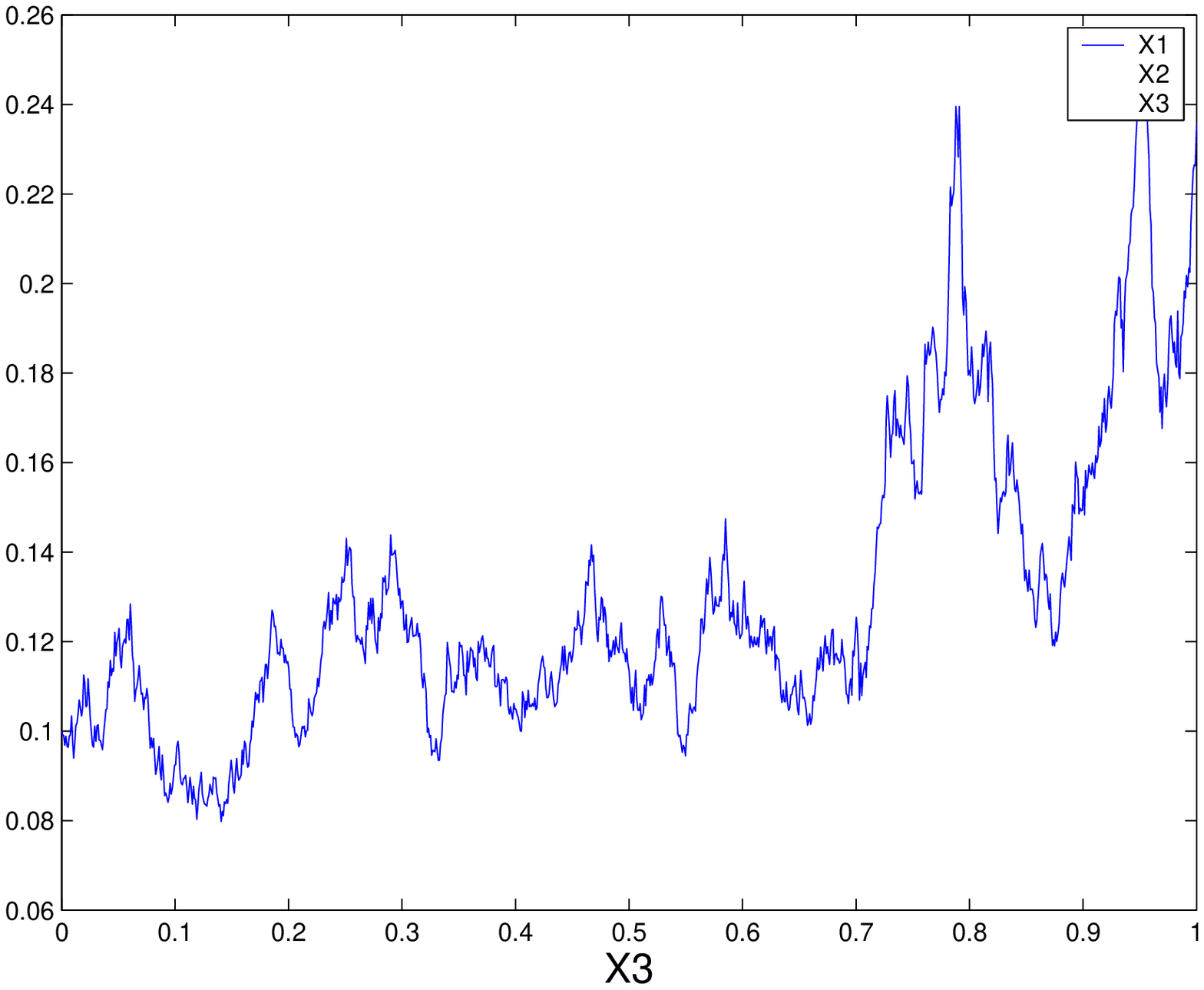}\\
  Pendulum with delay $\tau=1, \, H=0.5$ &  Euler top system with delay $\tau=1 \, H=0.5$ \\
\end{tabular}
\end{center}

For $K=0.3,$ and $\tau=1,$ the pendulum equation with delay
argument and with initial condition $\theta(0)=2,$ is represented
in the following figure. The Euler top system with delay argument
(\ref{1d}) is represented in the second figure, for the initial
conditions $x_1(0)=0.1, \, x_2(0)=0.05, \, x_3(0)=0.2.$
\begin{center}\begin{tabular}{cc}
  \includegraphics[height=2in]{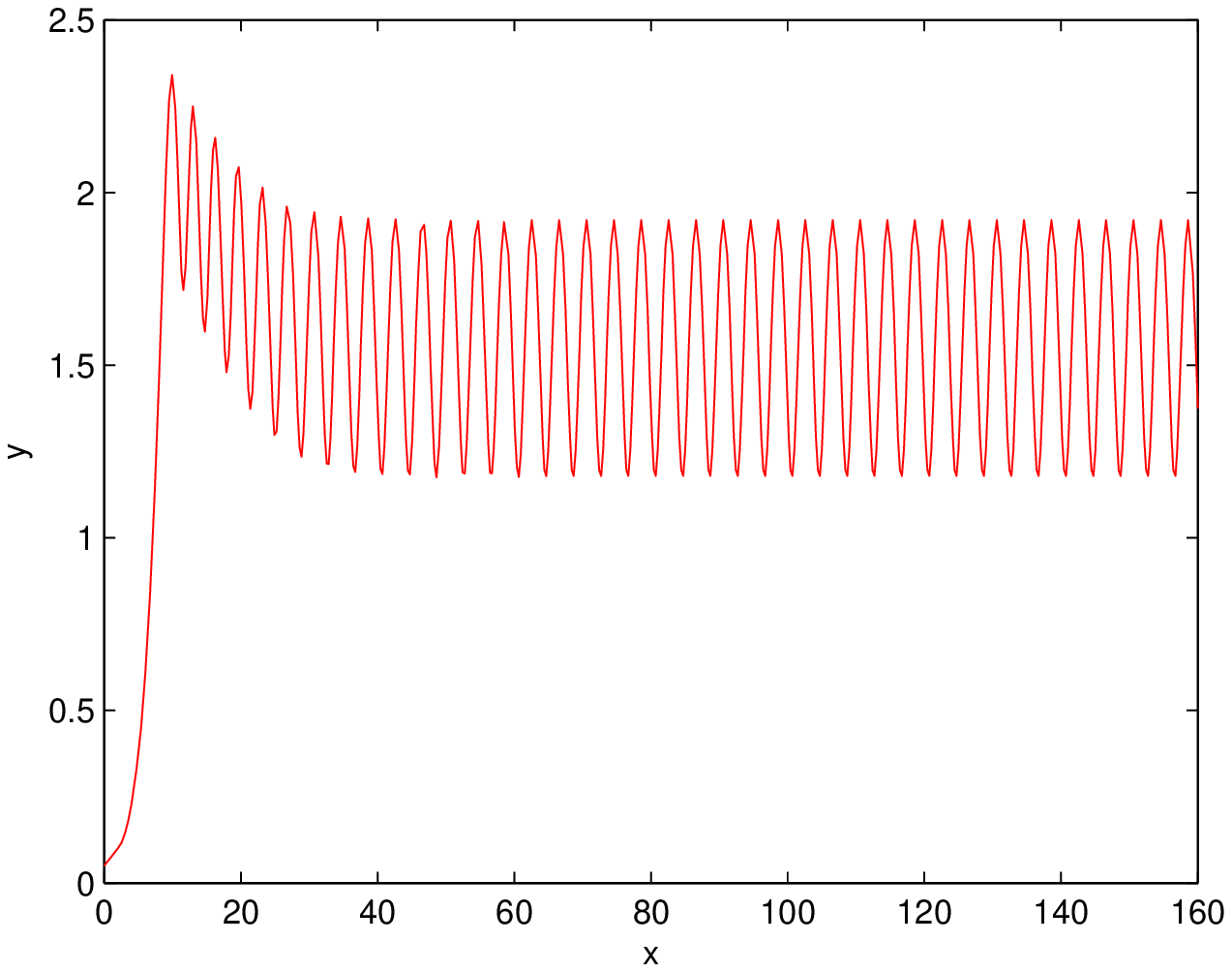}&
  \includegraphics[height=2in]{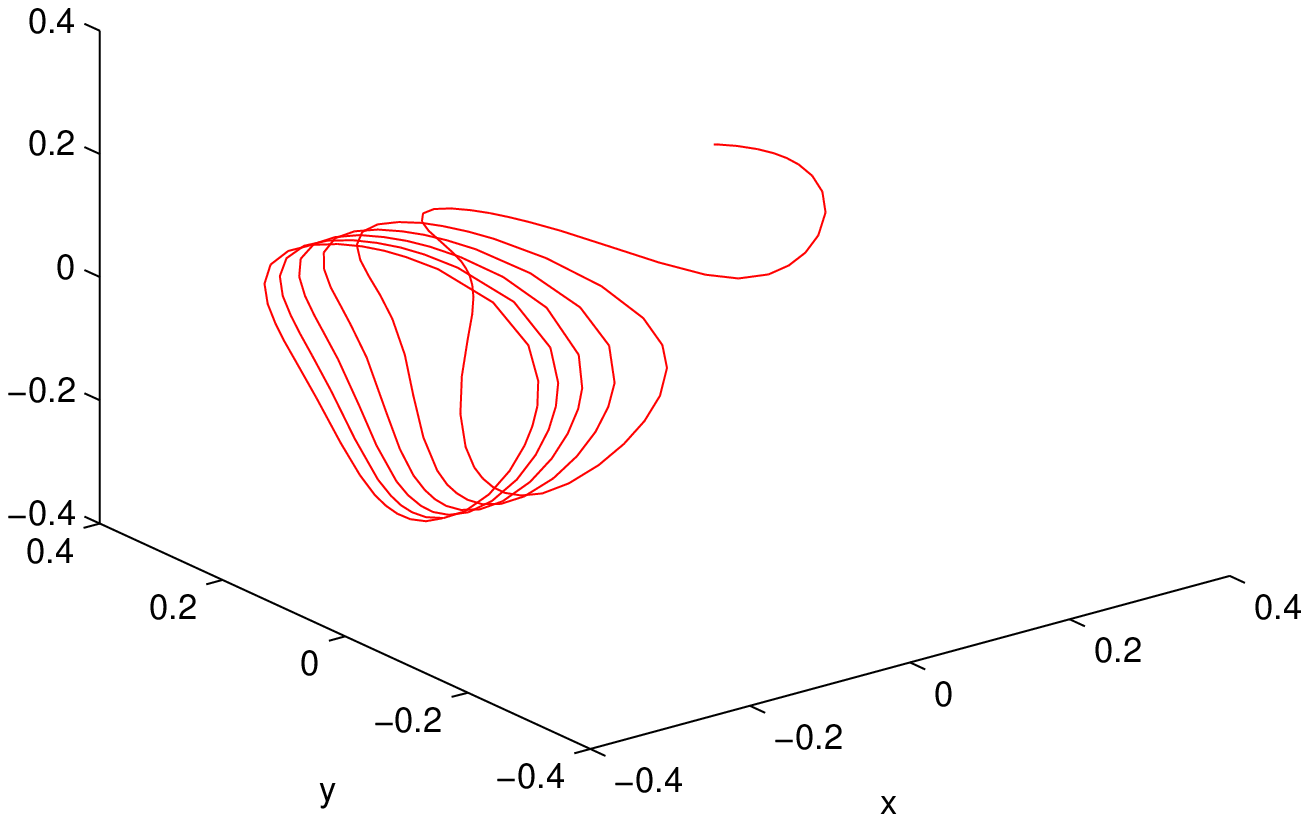}\\
Pendulum with delay $\tau=1, \, K=0.3$ &  Euler top system with delay $\tau=1 \, K=0.3$ \\
\end{tabular}
\end{center}

Using Caputo fractional derivative \cite{Di}, the following
propositions take place.

\begin{pr}
The Euler top system of fractional differential equations, given
by
\begin{equation}\label{1f}
\left\{%
\begin{array}{ll}
    \dot{x_{1}}(t)=x_{2}(t)x_{3}(t), &  \\
    \dot{x_{2}}(t)=-x_{1}(t)x_{3}(t), &  \\
    D^{\alpha}x_{3}(t)=x_{1}(t)x_{2}(t),&  \\
\end{array}%
\right.
\end{equation}with $\alpha\in(0,1),$ has the following properties
\begin{description}
    \item[a)] The function $H$ is a conservation law for
    \emph{(\ref{1f})};
    \item[b)] The solution of the system \emph{(\ref{1f})} on the
    constant level surface \emph{(\ref{c1})}, with $\theta(t)$ is
    the solution of the fractional equation
    \begin{equation}\label{pend1}
D^{\alpha+1}_t\theta(t)+2H\sin\theta(t)=0,
\end{equation} and reciprocal.
\end{description}
\end{pr}\hfill $\Box$

\begin{pr}
The Euler top system of fractional differential equations, given
by
\begin{equation}\label{2f}
\left\{%
\begin{array}{ll}
    D^{\alpha}x_{1}(t)=x_{2}(t)x_{3}(t), &  \\
    \dot{x}_{2}(t)=-x_{1}(t)x_{3}(t), &  \\
    \dot{x}_{3}(t)=x_{1}(t)x_{2}(t),&  \\
\end{array}%
\right.
\end{equation}with $\alpha\in(0,1),$ has the following properties
\begin{description}
    \item[a)] The function $H$ is a conservation law for
    \emph{(\ref{2f})};
    \item[b)] The solution of the system \emph{(\ref{2f})} on the
    constant level surface \emph{(\ref{c2})}, with $\theta(t)$ is
    the solution of the fractional equation
    \begin{equation}\label{pend2}
D^{\alpha+1}_t\theta(t)+2K\sin\theta(t)=0,
\end{equation} and reciprocal.
\end{description}
\end{pr}\hfill $\Box$

By using the Adams-Moulton method for integration, for the initial
condition $\theta(0)=-3.1,$ the solution of the fractional
differential equation (\ref{pend1}) is represented in the
fo\-llowing gra\-phics for $\alpha =0.8,$ respectively for
$\alpha=1.$

\begin{center}\begin{tabular}{ccc}
  \includegraphics[height=2in]{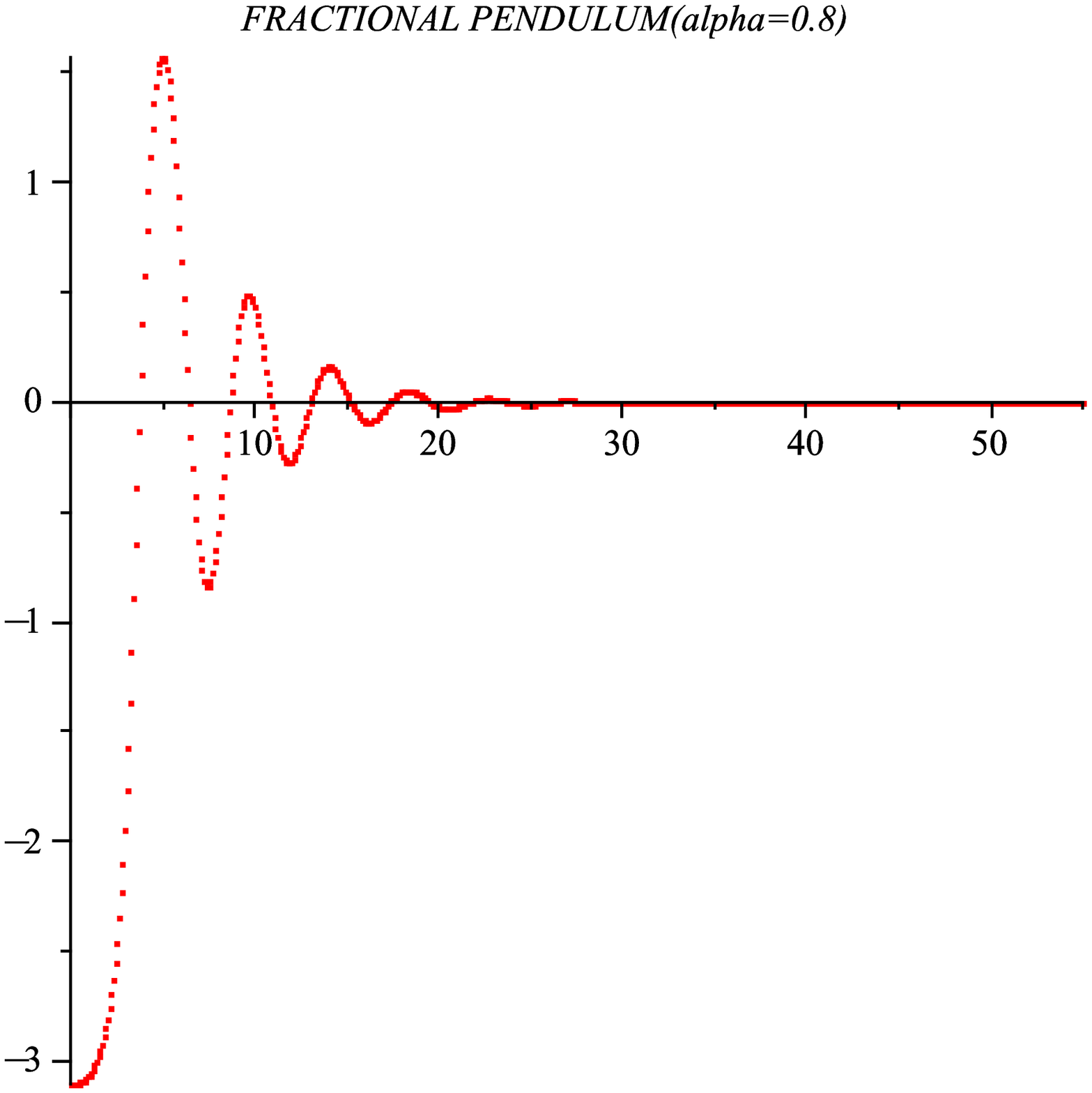}
  \includegraphics[height=2in]{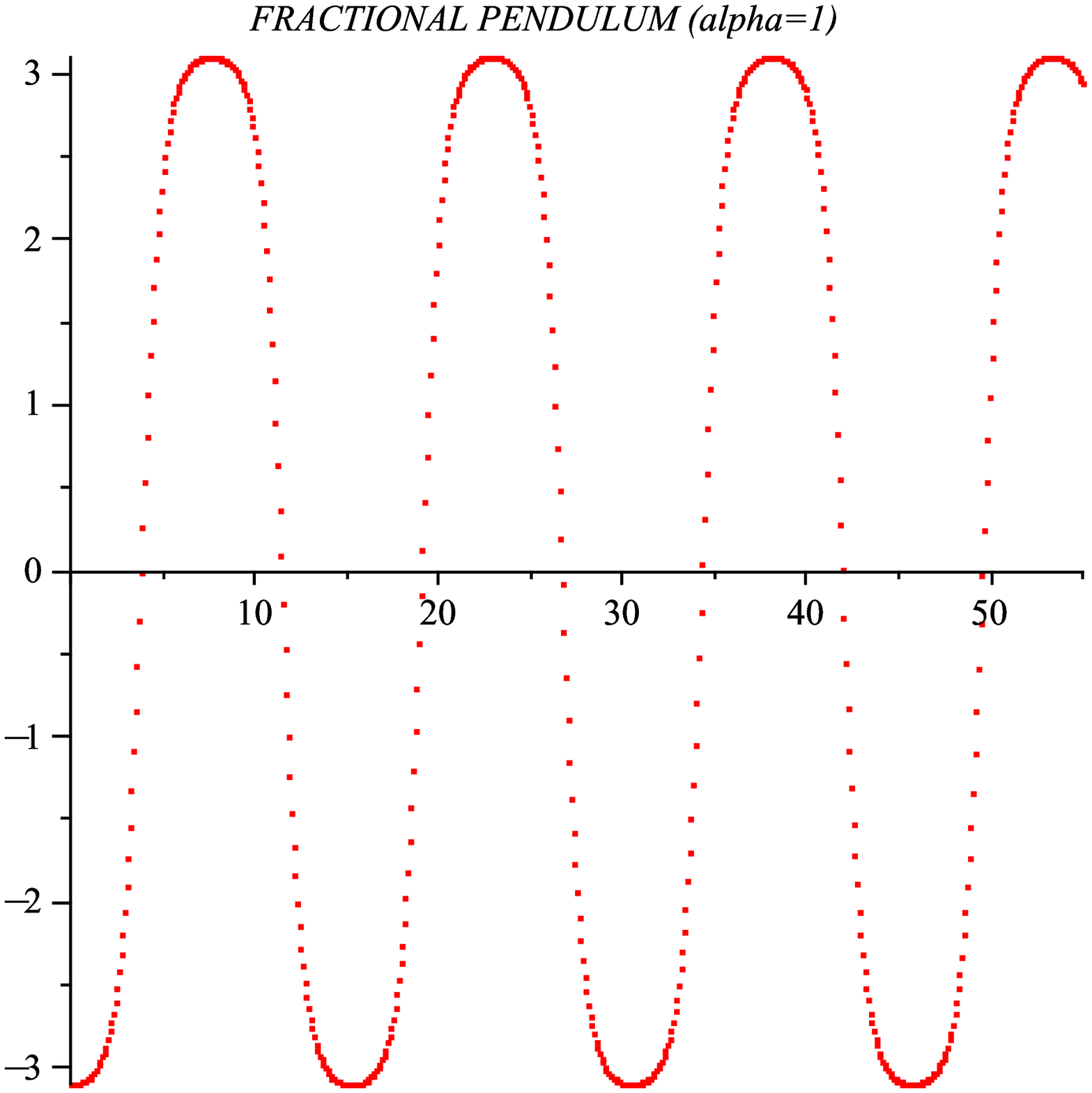}
\end{tabular}
\end{center}

It can be observed that the pendulum solution is asymptotically
stable for  $0<\alpha<1$ and it is oscillatory for $\alpha=1.$

The solution for the system of fractional differential equations
(\ref{1f}), respectively for (\ref{2f}), is represented in the
above graphics, for the initial conditions  $x_1(0)=0.1, \,
x_2(0)=0.1$ and $x_3(0)=0.3.$ The cases of  $\alpha=0.8$ and
$\alpha=1$ are illustrated.

\begin{center}\begin{tabular}{ccc}
  \includegraphics[height=2in]{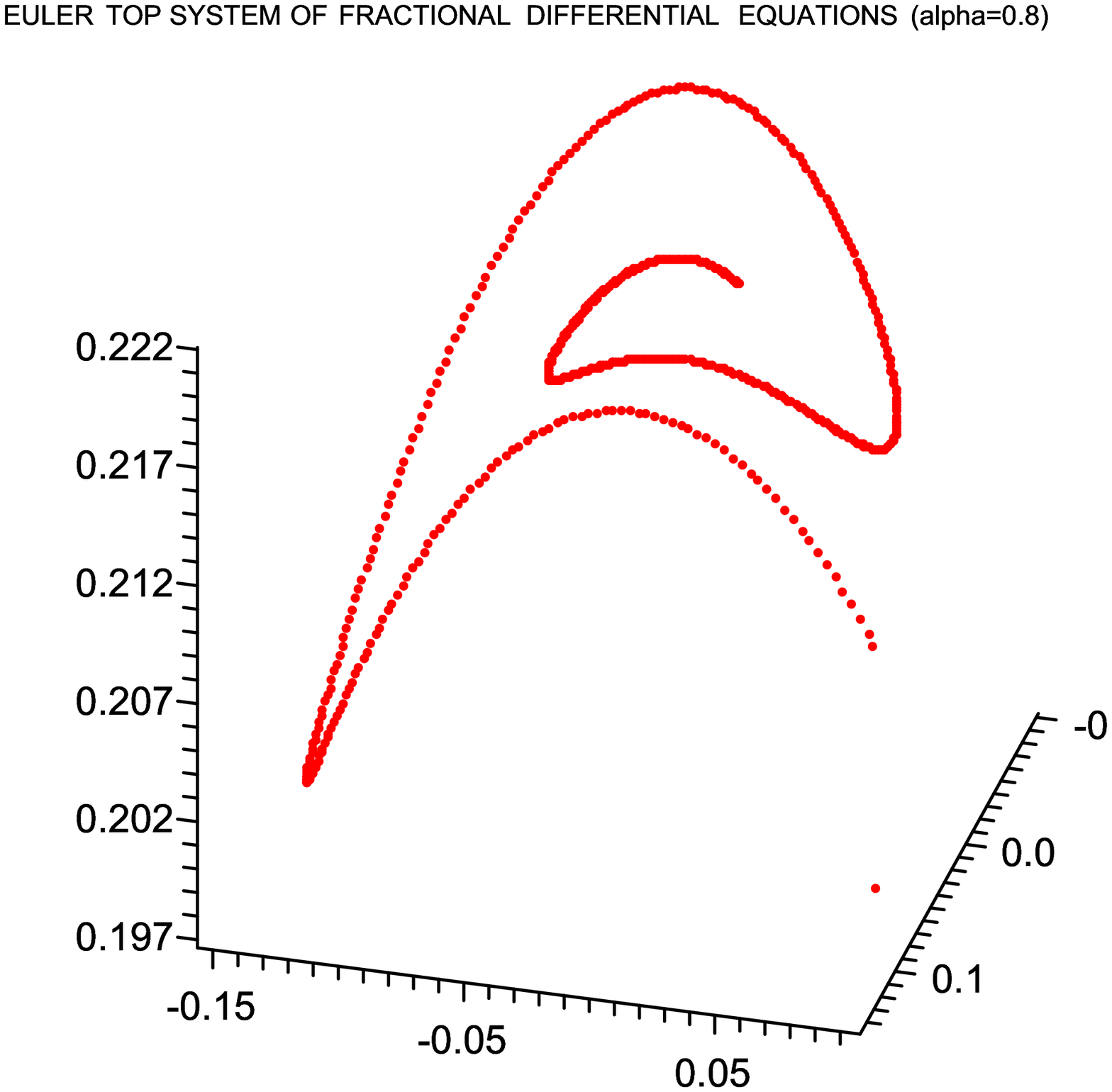}
  \includegraphics[height=2in]{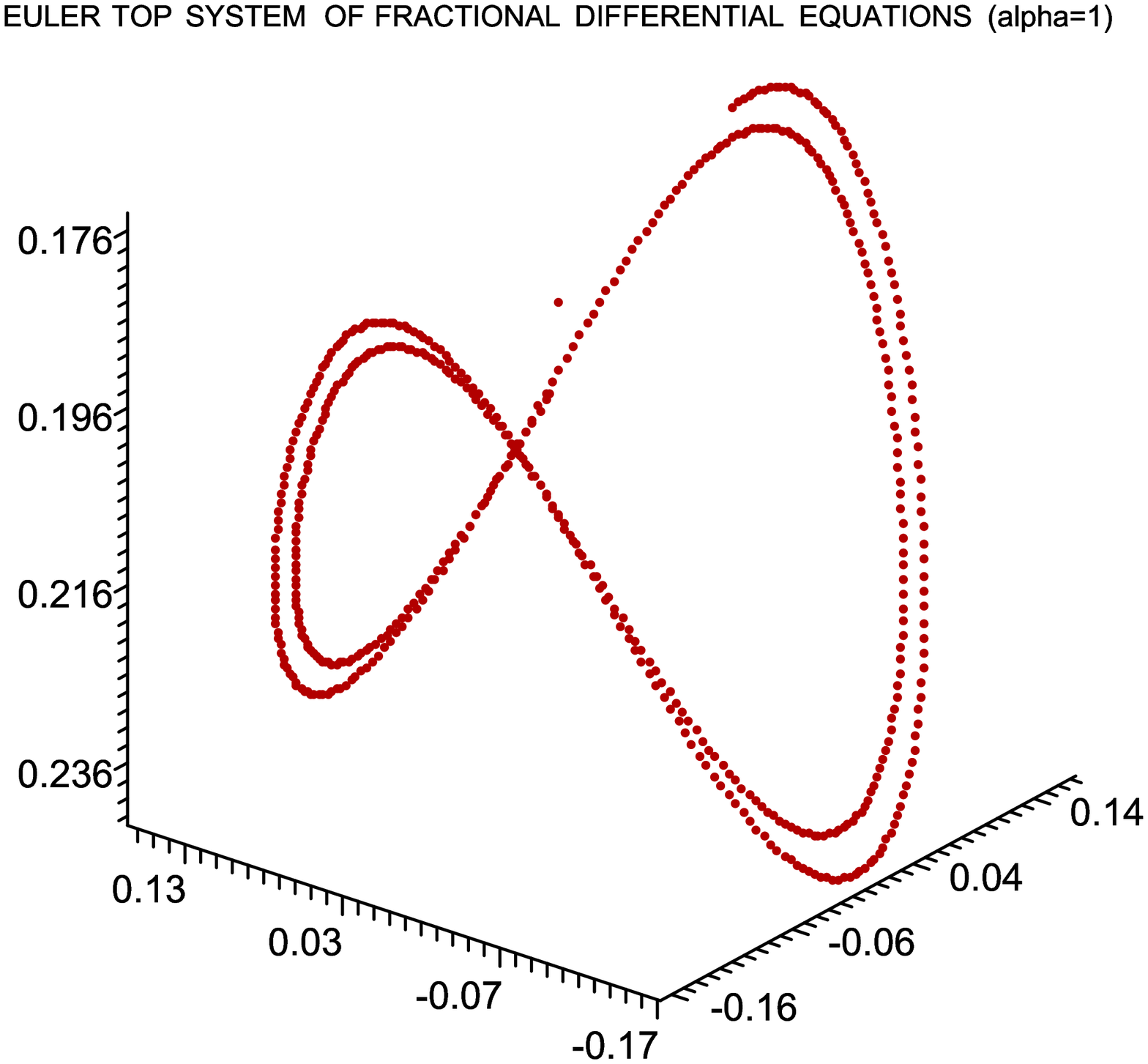}
\end{tabular}
\end{center}

\begin{center}\begin{tabular}{ccc}
  \includegraphics[height=2in]{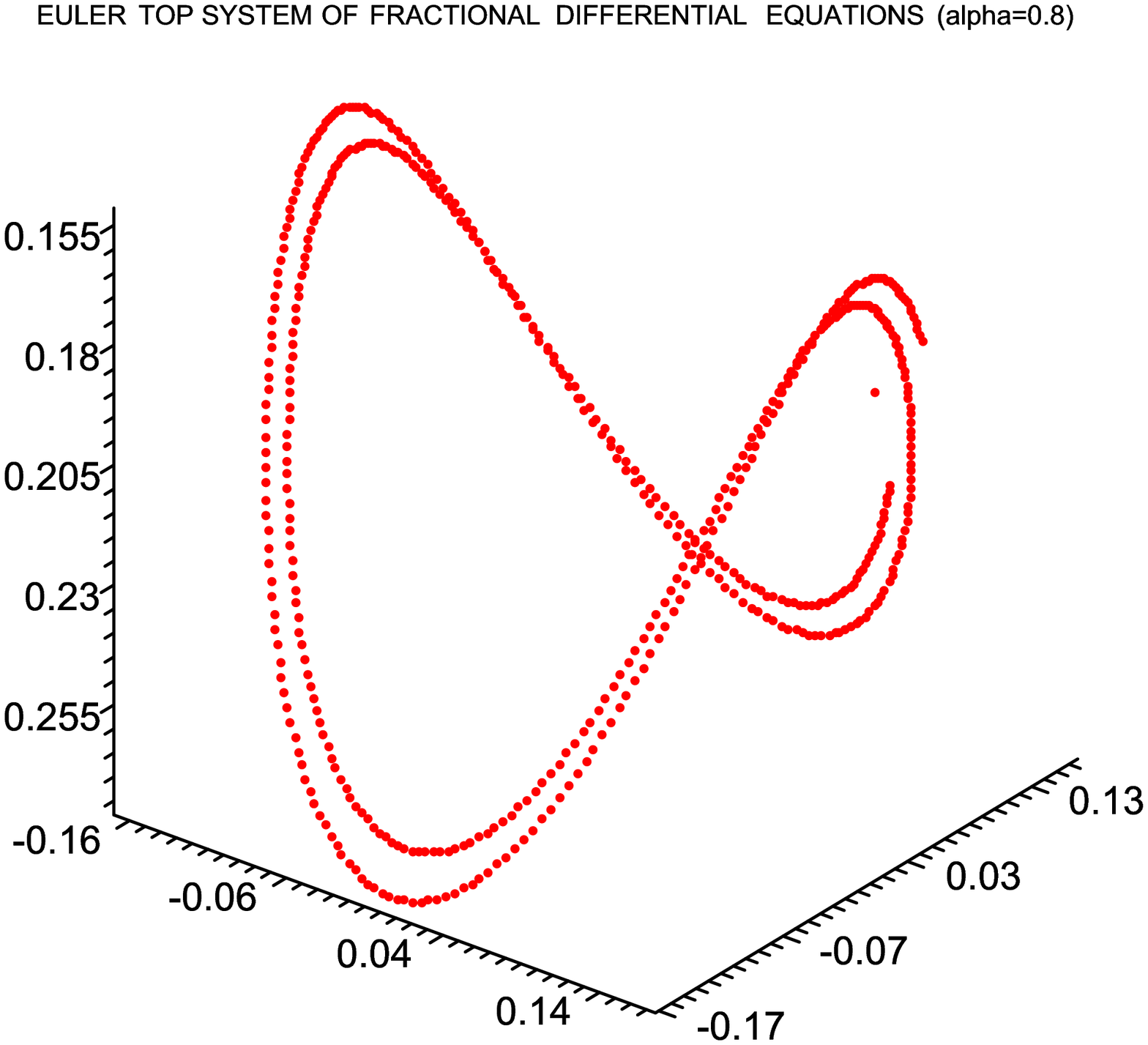}
  \includegraphics[height=2in]{SISEULERTOP1.eps}
\end{tabular}
\end{center}

\section{Stochastic Euler top system and stochastic pendulum}

A Wiener process describes rapidly  fluctuating random phenomena.
Stochastic di\-fferential equations (SDE) are stochastic integral
equations and are written symbo\-lically in a differential form.
We will consider such a Wiener process of the form
\begin{equation}\label{s1}
dx(t)=f(x(t))dt+g(x(t))dW(t),
\end{equation} where $f$ is the slowly varying continuous
component called drift coefficient and $g$ is the rapidly varying
continuous component called diffusion coefficient.
 The integral representation is of the form
\begin{equation}
x(t)=x(t_0)+\int_{t_0}^tf(x(s))ds+\int_{t_0}^t g(x(s))dW(s),
\end{equation}
where $W(t)$is a Wiener process, a Gaussian process with $W(0)=0$
and $N(0,t)-$distributed $W(t)$ for each $t\geq 0,$ so
$$\mathbb{E}(W(t))=0, \, \mathbb{E}((W(t))^2)=t.$$ The first
integral is a Riemann-Stieltjes integral and the second one is a
stochastic integral. The most studied interpretation of the
stochastic integral are those of It\^{o} and Stratonovich. The
choice of interpretation depends on the type of analysis required
for solution \cite{High}. It\^{o} stochastic calculus is closely
related to diffusion processes and martingale theory \cite{High}.
The solution of (\ref{s1}) is a diffusion process with transition
probability $p=u(x(t)),$ satisfying the Fokker-Planck equation
\begin{equation}\label{f-p}
\frac{\partial}{\partial t}u(x(t))=-\frac{\partial}{\partial
x(t)}f(x(t)) u(x(t))+\frac{1}{2}\frac{\partial^2}{\partial
(x(t))^2}[(g(x(t)) g^T(x(t))u(x(t))].
\end{equation} Equations (\ref{s1}) and (\ref{f-p}) contain the
same statistical information from a one-particle process point of
view (but not if we think the It\^{o} equation as describing a
random dynamical system) \cite{sim}.

An It\^{o} SDE is written in the form (\ref{s1}) and a
Stratonovich SDE is written symbolically in the form
\begin{equation}\label{s2}
dx(t)=f(x(t))dt+g(x(t))\circ dW(t),
\end{equation}and in the integral form as
\begin{equation}
x(t)=x(t_0)+\int_{t_0}^tf(x(s))ds+\int_{t_0}^t g(x(s))\circ dW(s).
\end{equation}

It is possible to switch between these two approaches, in the
sense that the It\^{o} SDE (\ref{s1}) has the same solution as
Stratonovich SDE
\begin{equation}\label{s3}
d x(t)=\underline{f}(x(t))dt+g(x(t))\circ dW(t),
\end{equation}with modified drift coefficient
$$\underline{f}(x(t))=f(x(t))-\frac{1}{2}g(x(t))\frac{\partial g}{\partial x(t)}(x(t)).$$

If $W(1), ..., W(d)$ are $d$ independent Wiener processes, and
$x(t)=(x_1(t),...,x_n(t))$ then the multi-Wiener process case can
be written in the form
\begin{equation}
dx_i(t)=f^i(x(t))dt+\sum_{j=1}^d g_{ij}(x(t))dW(j),
\end{equation}with $g(x(t))$ a $n\times d$ matrix and $dW$ a $d\times
1$ matrix.

In Stratonovich case, the stochastic system of differential
equations with a multi-Wiener process, can be written in the
following manner
\begin{equation}\label{strat}
dx(t)=\underline{f}(x(t))dt+\sum_{j=1}^d g_j(x(t))\circ dW(j),
\end{equation}where
$$\underline{f}(x(t))=f(x(t))-\frac{1}{2}\sum_{k=1}^{n}\sum_{j=1}^d g_{k,j}(x(t))
\frac{\partial g_j}{\partial x_k(t)}.$$

The Euler top system of stochastic differential equations can be
represented in the following form,
\begin{equation}\label{s4}
\left\{%
\begin{array}{ll}
    dx_1(t)=x_2(t)x_3(t)dt+x_1(t)dW^1(t) &  \\
    dx_2(t)=-x_{1}(t)x_{3}(t)dt, &  \\
    dx_3(t)=x_{1}(t)x_{2}(t)dt+dW^3(t),&  \\
\end{array}%
\right.
\end{equation}with the Wiener process $W(t)=(W^1(t),0,W^2(t)),$ the
drift coefficients $f^1(x(t))=x_2(t)x_3(t),
f^2(x(t))=-x_1(t)x_3(t), \, f^3(x(t))=x_1(t)x_2(t),$
$x(t)=(x_1(t),x_2(t),x_3(t))^T, \\
f(x(t))=(f^1(x(t)),f^2(x(t)),f^3(x(t)))^T$ and the diffusion
coefficient vectors
$$g^1(x(t))= \left( \begin{array}{ccc}
x_1(t)  \\
0 \\
0
\end{array} \right), \, g^2(x(t))= \left( \begin{array}{ccc}
0  \\
0 \\
0
\end{array} \right), \, g^3(x(t))= \left( \begin{array}{ccc}
0  \\
0 \\
1
\end{array} \right).$$

The corresponding (It\^{o}) Fokker-Planck equation for the
probability density $p=u(x(t))$ reads
\begin{eqnarray*}
\frac{\partial}{\partial t}u(x(t))&=&-\frac{\partial}{\partial
x_1(t)}[x_2(t)x_3(t)u(x(t))]+\frac{\partial}{\partial
x_2(t)}[x_1(t)x_3(t) u(x(t))]-\frac{\partial}{\partial
x_3(t)}[x_1(t)x_2(t) u(x(t))]\\
&+&\frac{1}{2}\frac{\partial^2}{\partial
(x_1(t))^2}[(x_1(t))^2u(x(t))]+\frac{1}{2}\frac{\partial^2}{\partial
(x_3(t))^2}u(x(t)).
\end{eqnarray*}

In the Stratonovich case, stochastic system (\ref{s4}) can be
written using relation (\ref{strat}) in the following manner
\begin{equation}\label{s4s}
\left\{%
\begin{array}{ll}
    dx_1(t)=(x_2(t)x_3(t)-\frac{1}{2}x_2(t))dt+x_1(t)dW^1(t) &  \\
    dx_2(t)=-x_{1}(t)x_{3}(t)dt, &  \\
    dx_3(t)=x_{1}(t)x_{2}(t)dt+dW^3(t).&  \\
\end{array}%
\right.
\end{equation}

The stochastic system (\ref{s4}), respectively (\ref{s4s}),  is
implemented in Matlab, using Milstein scheme, for initial
conditions $x_1(1)=0.1, \, x_2(1)=0.1, \, x_3(1)=0.1$ and orbits
are represented in the following figures.

\begin{center}\begin{tabular}{ccc}
  \includegraphics[height=2in]{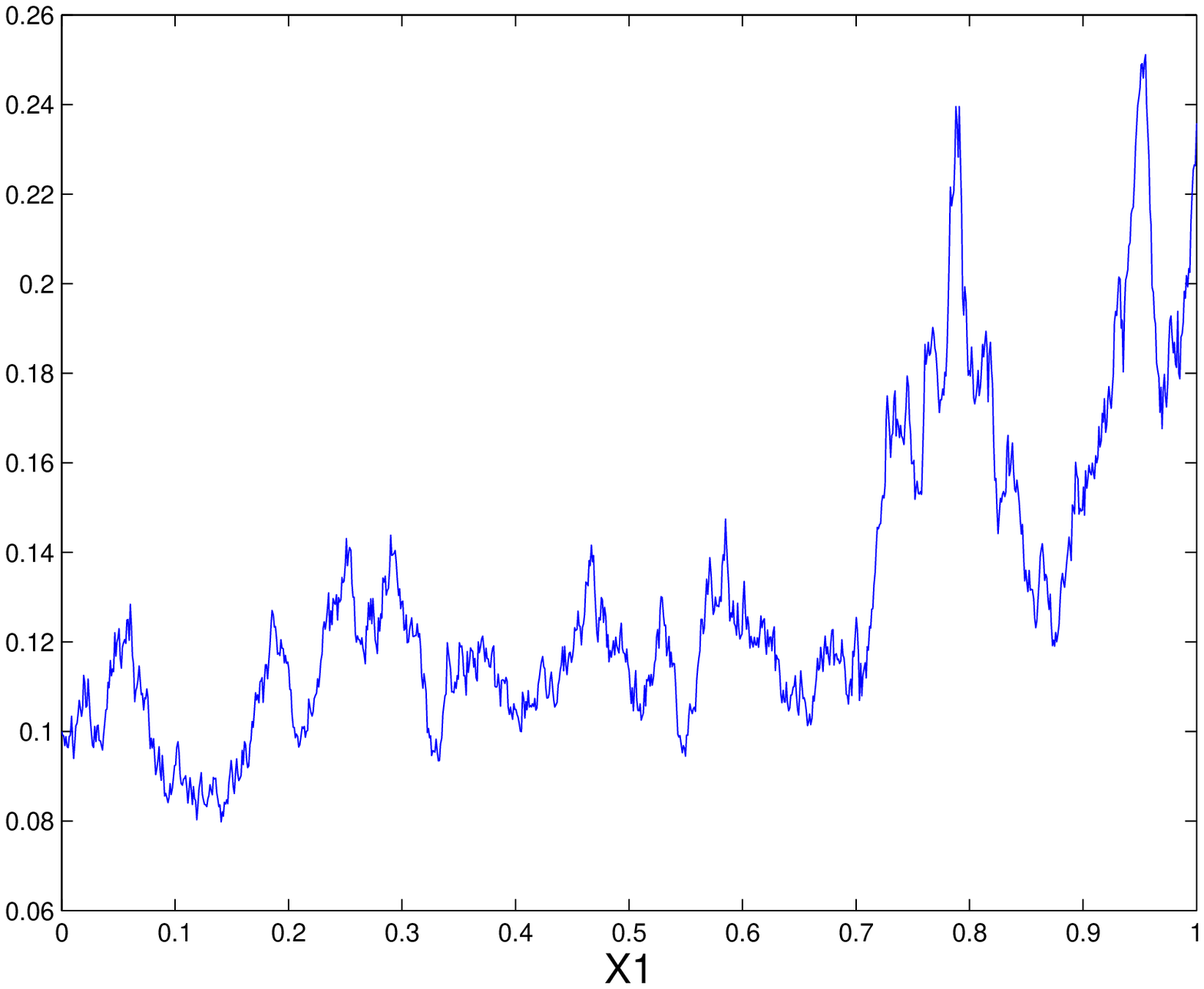}
  \includegraphics[height=2in]{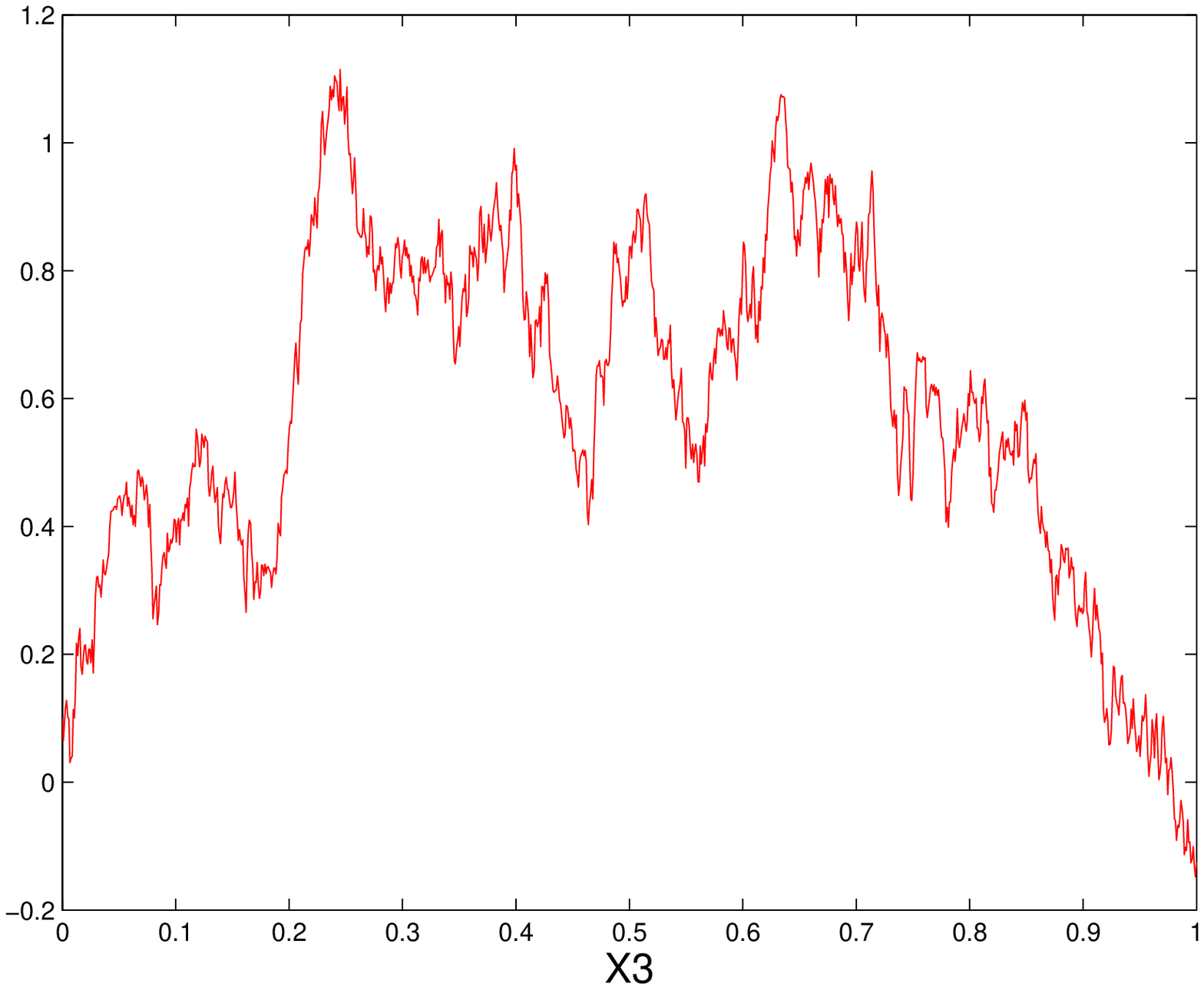}
\end{tabular}
\end{center}

\begin{center}\begin{tabular}{ccc}
  \includegraphics[height=2in]{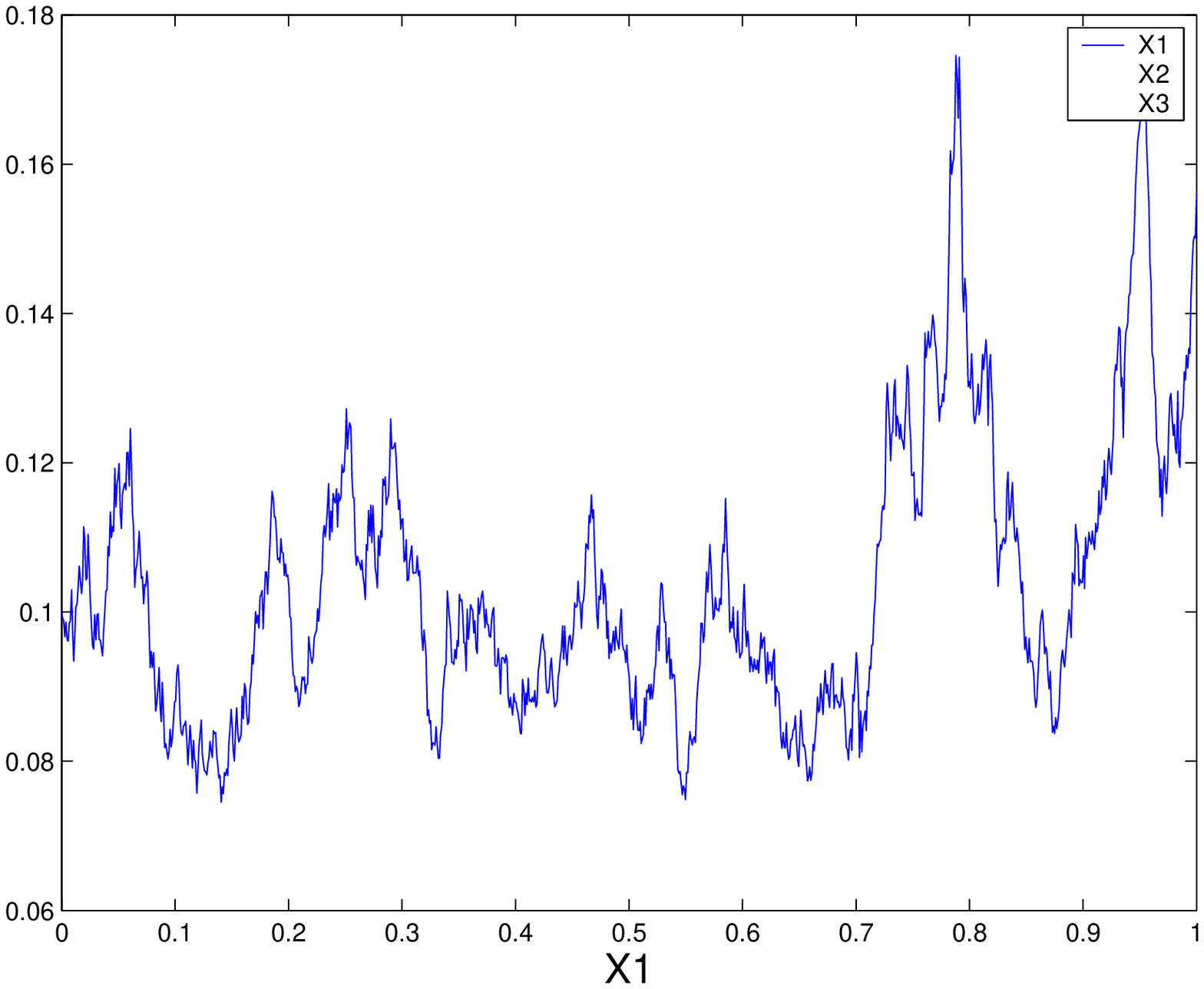}
  \includegraphics[height=2in]{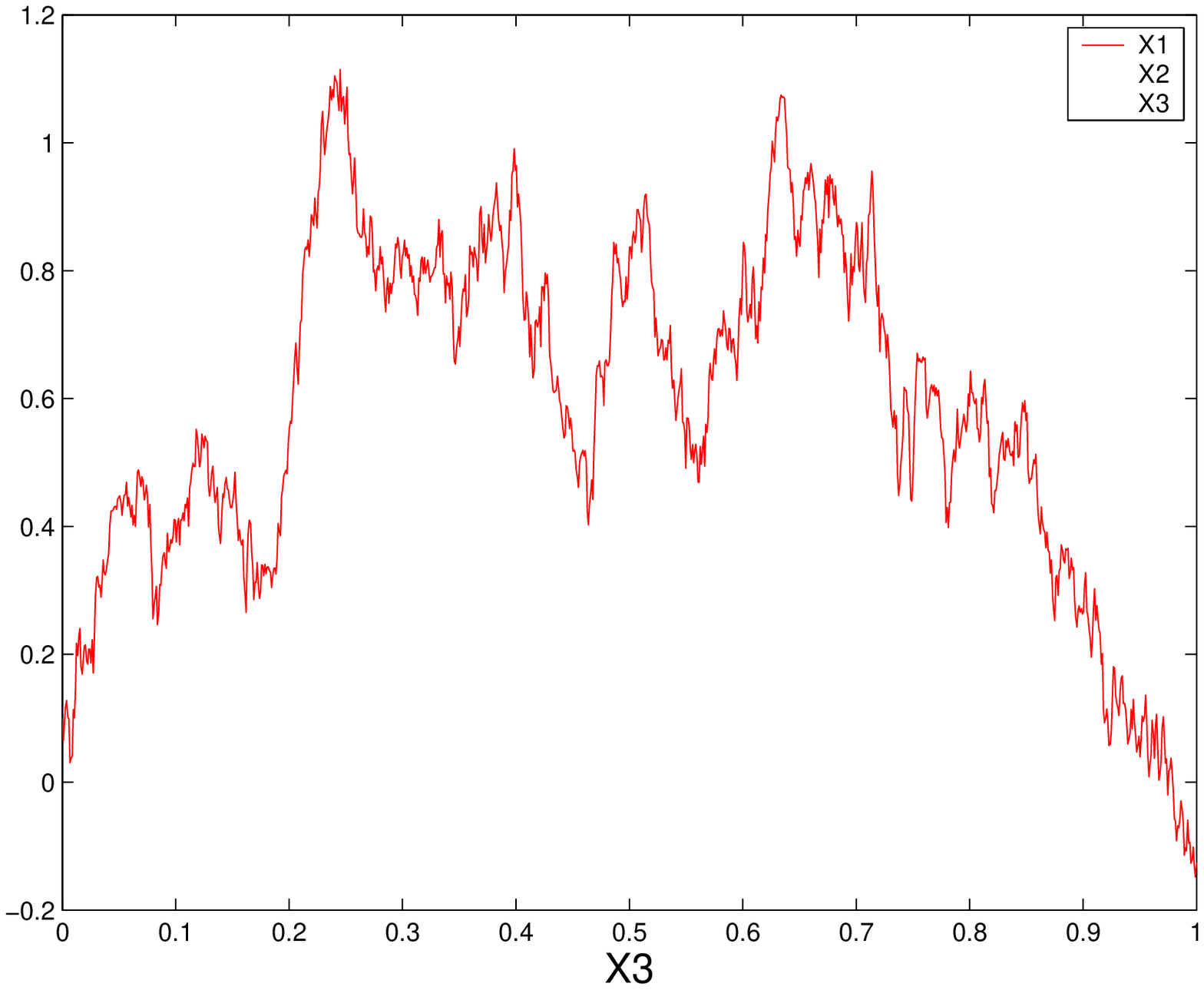}
\end{tabular}
\end{center}

If the SDE of Euler top system has the form
\begin{equation}\label{s51}
\left\{%
\begin{array}{ll}
    dx_1(t)=x_2(t)x_3(t)dt+\sqrt{x_1(t)}dW^1(t) &  \\
    dx_2(t)=-x_1(t)x_3(t)dt+\sqrt{x_2(t)}dW^2(t), &  \\
    dx_3(t)=x_1(t)x_2(t)dt+\sqrt{x_3(t)}dW^3(t),&  \\
\end{array}%
\right.
\end{equation}
then drift coefficients are $$f^1((t))=x_2(t)x_2(t), \,
f^2(x(t))=-x_1(t)x_3(t), \, f^3(x(t))=x_1(t)x_2(t),$$ with
$x(t)=(x_1(t),x_2(t),x_3(t))^T, \,
f(x(t))=(f^1(x(t)),f^2(x(t)),f^3(x(t)))^T,$ and the diffusion
coefficient vectors
$$g^1(x(t))= \left( \begin{array}{ccc}
\sqrt{x_1(t)}  \\
0 \\
0
\end{array} \right), \, g^2(x(t))= \left( \begin{array}{ccc}
0  \\
\sqrt{x_2(t)} \\
0
\end{array} \right), \, g^3(x(t))= \left( \begin{array}{ccc}
0  \\
0 \\
\sqrt{x_3(t)}
\end{array} \right),$$ then
the associated (It\^{o}) Fokker-Planck equation for the
probability density $p=u(x(t))$ is
\begin{eqnarray*}
\frac{\partial}{\partial t}u(x(t))&=&-\frac{\partial}{\partial
x_1(t)}[x_2(t)x_3(t)u(x(t))]+\frac{\partial}{\partial
x_2(t)}[x_1(t)x_3(t) u(x(t))]-\frac{\partial}{\partial
x_3(t)}[x_1(t)x_2(t) u(x(t))]\\
&+&\frac{1}{2}\frac{\partial^2}{\partial
(x_1(t))^2}[x_1(t)u(x(t))]+\frac{1}{2}\frac{\partial^2}{\partial
(x_2(t))^2}[x_2(t)u(x(t))]+\frac{1}{2}\frac{\partial^2}{\partial
(x_3(t))^2}[x_3(t)u(x(t))].
\end{eqnarray*}

The Stratonovich stochastic Euler top system is written in the
following way
\begin{equation}\label{s51s}
\left\{%
\begin{array}{ll}
    dx_1(t)=(x_2(t)x_2(t)-\frac{1}{4})dt+\sqrt{x_1(t)}dW^1(t) &  \\
    dx_2(t)=-(x_1(t)x_3(t)+\frac{1}{4})dt+\sqrt{x_2(t)}dW^2(t), &  \\
    dx_3(t)=(x_1(t)x_2(t)dt+\frac{1}{4})+\sqrt{x_3(t)}dW^3(t),&  \\
\end{array}%
\right.
\end{equation}

Stochastic system (\ref{s51}), respectively (\ref{s51s}), can be
implemented using stochastic Euler method which represents a
square-root model. For initial values $x_1(1)=1, \, x_2(1)=0.8, \,
x_3(1)=0.2,$ orbits are represented in the following figures.

\begin{center}\begin{tabular}{ccc}
  \includegraphics[height=1.5in]{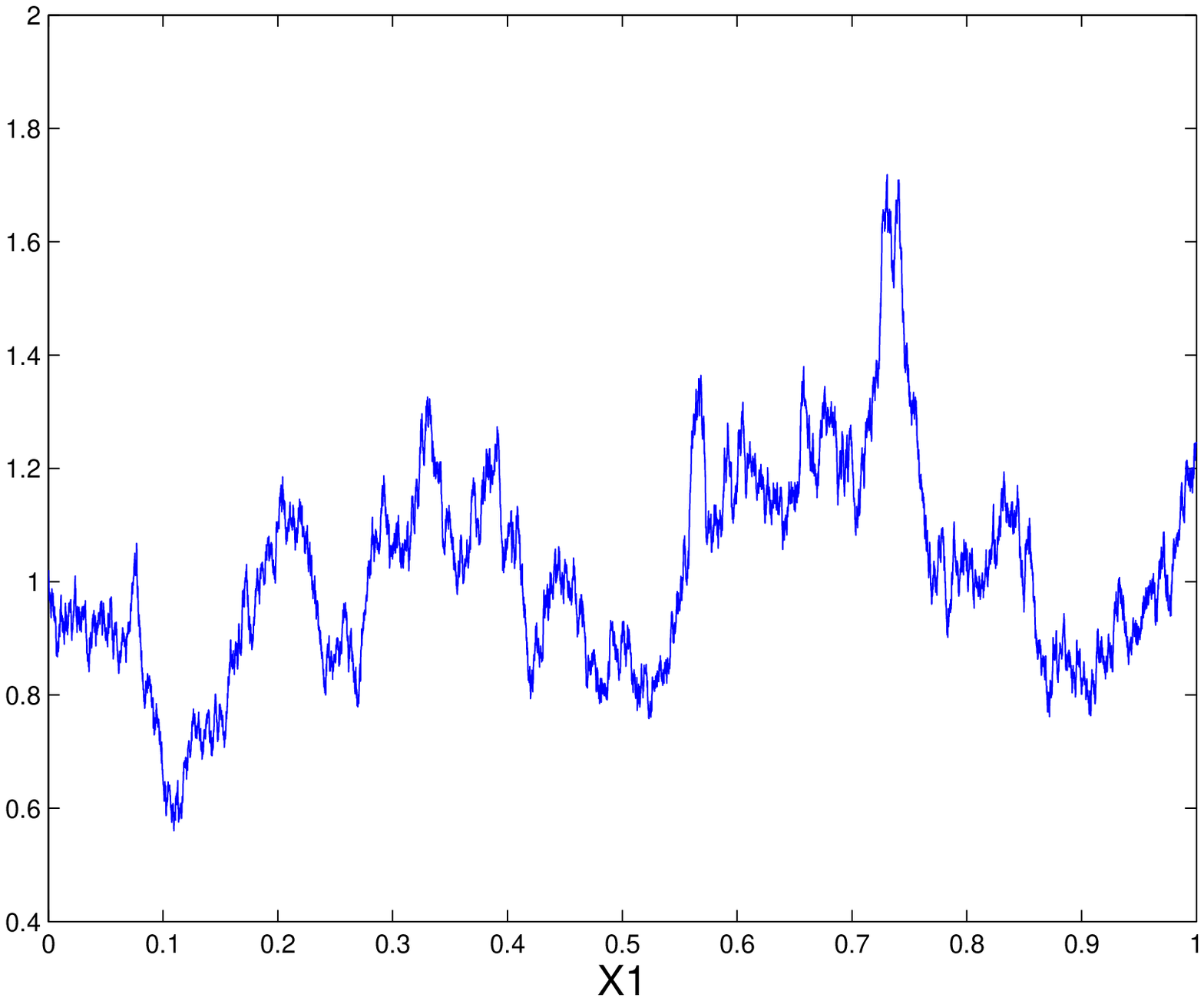}
  \includegraphics[height=1.5in]{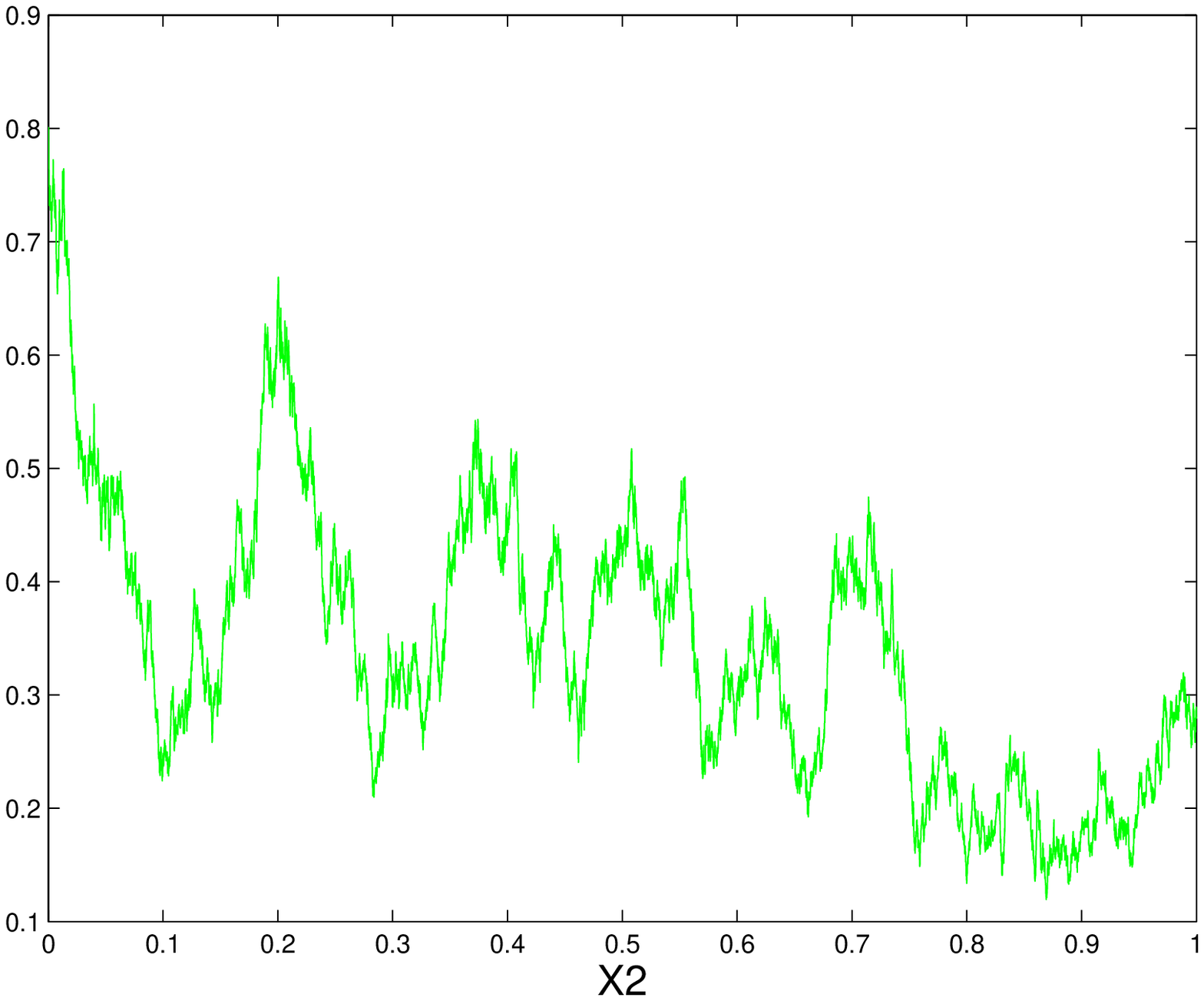}
  \includegraphics[height=1.5in]{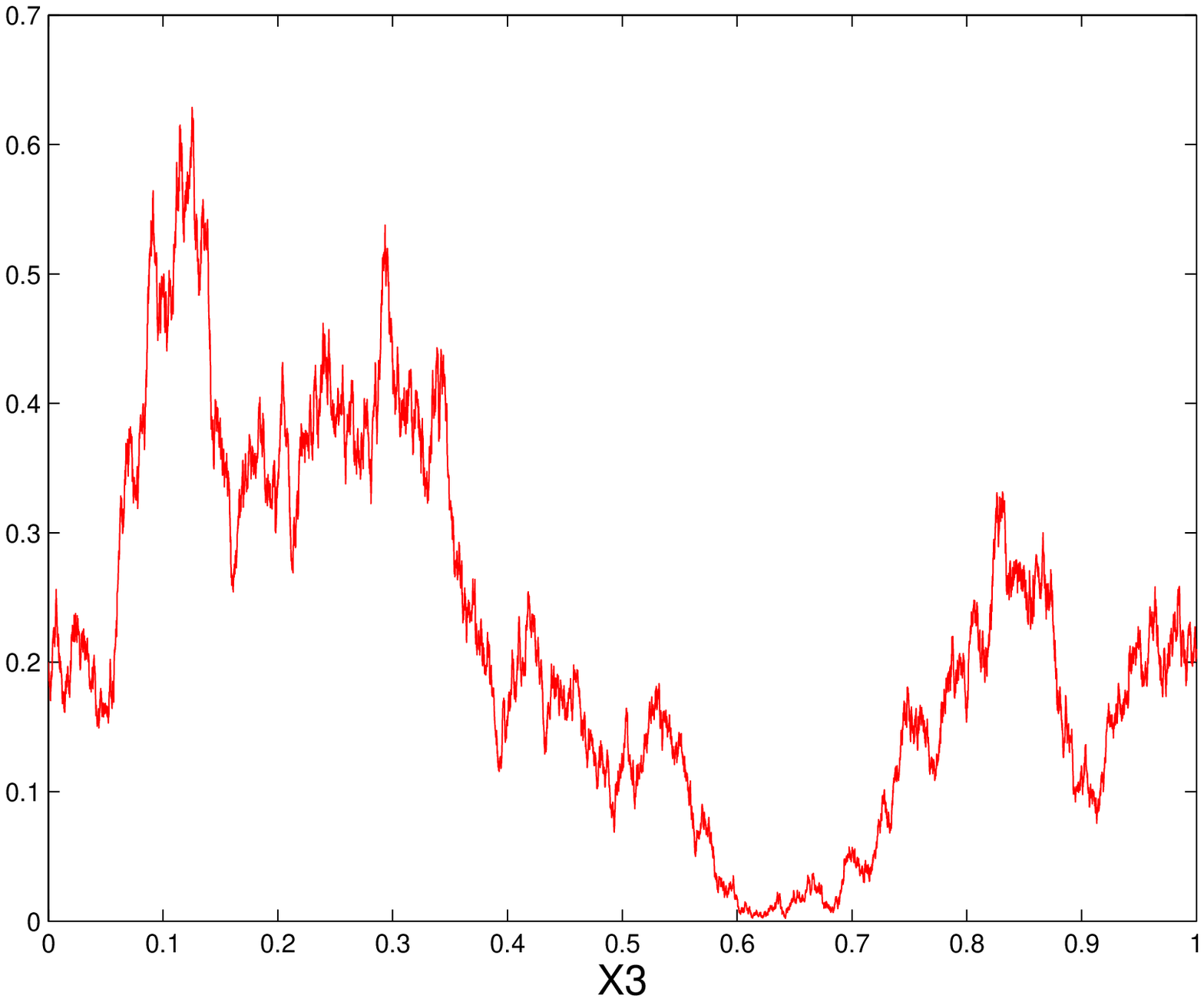}
\end{tabular}
\end{center}

\begin{center}\begin{tabular}{ccc}
  \includegraphics[height=1.5in]{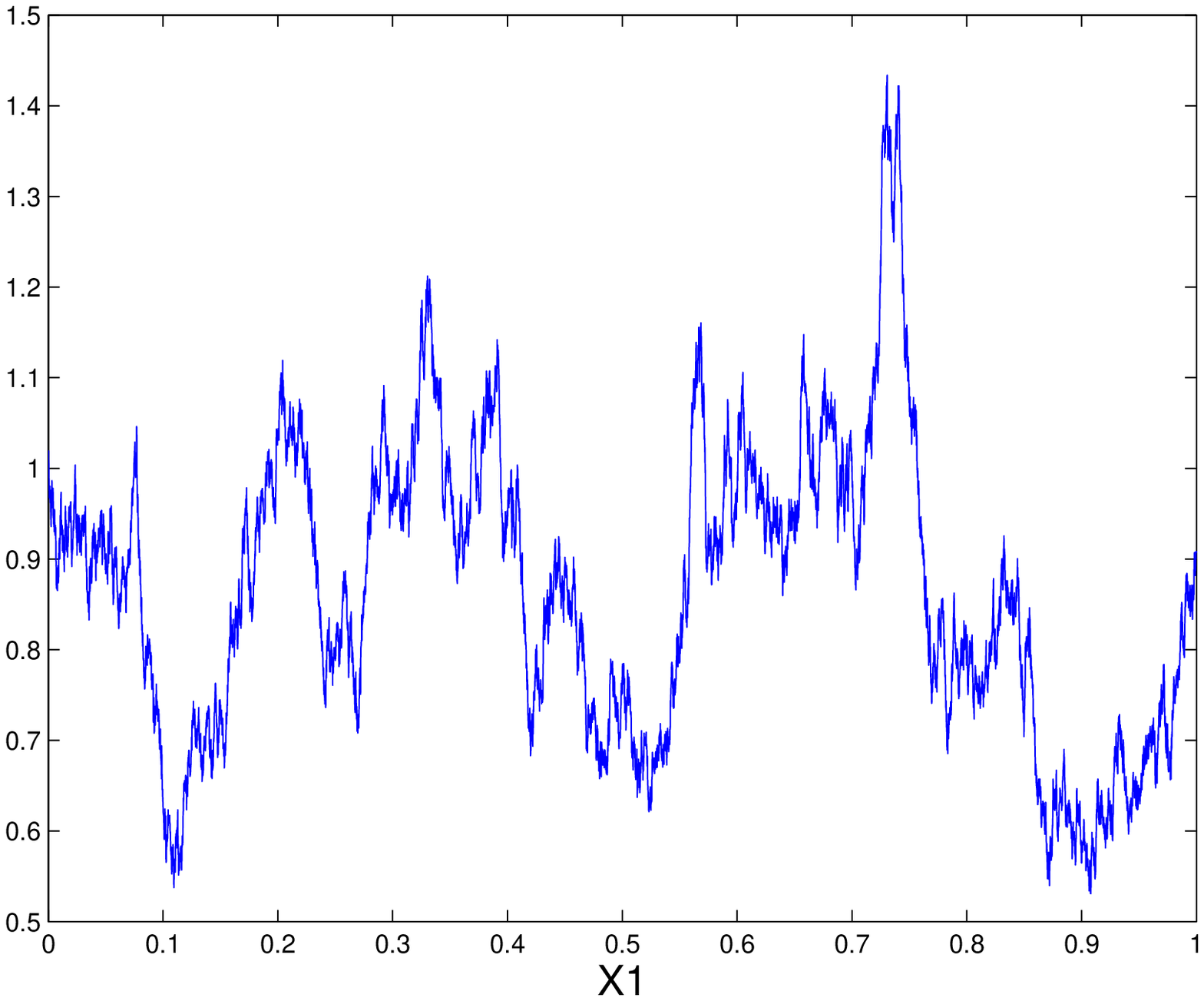}
  \includegraphics[height=1.5in]{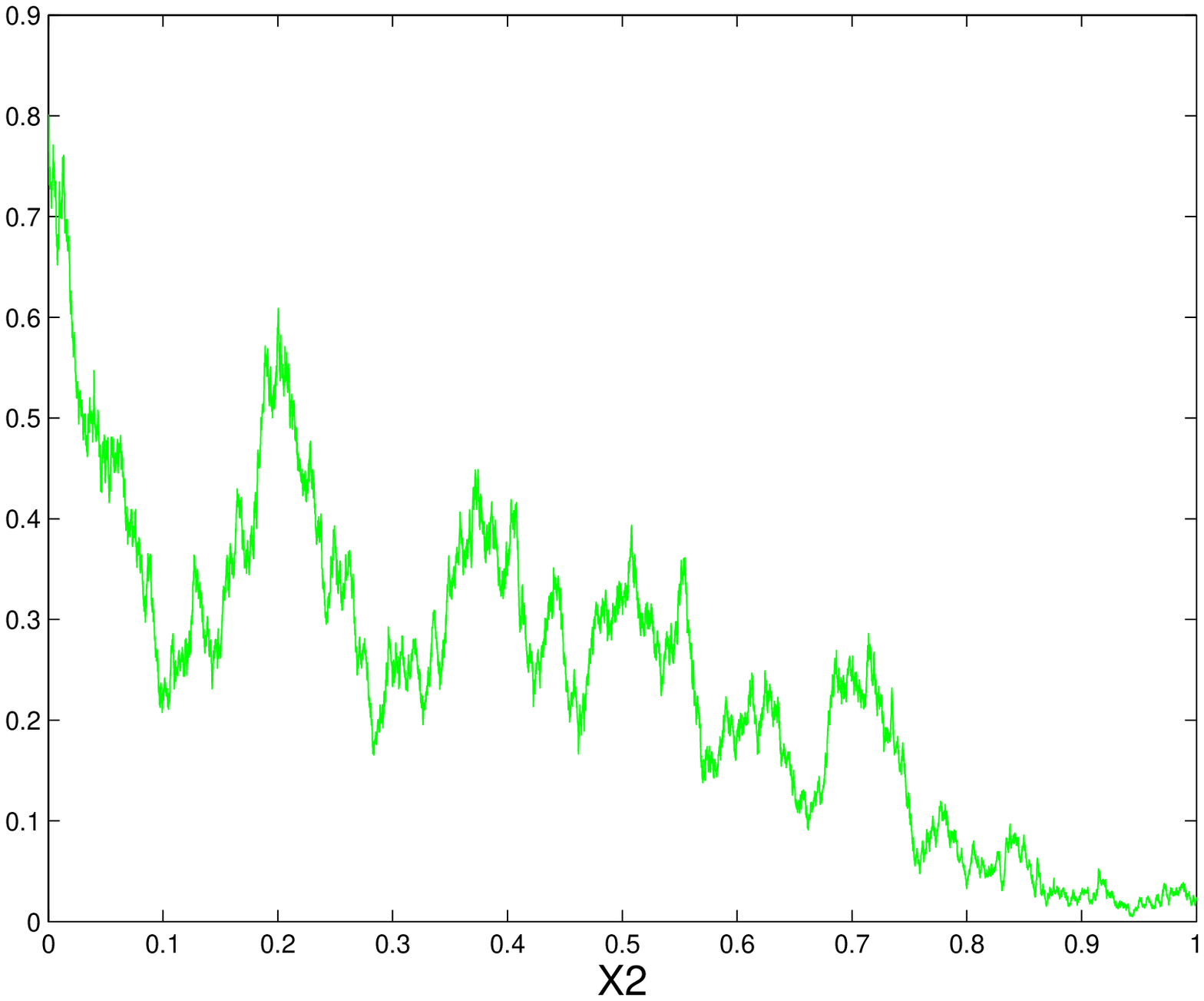}
  \includegraphics[height=1.5in]{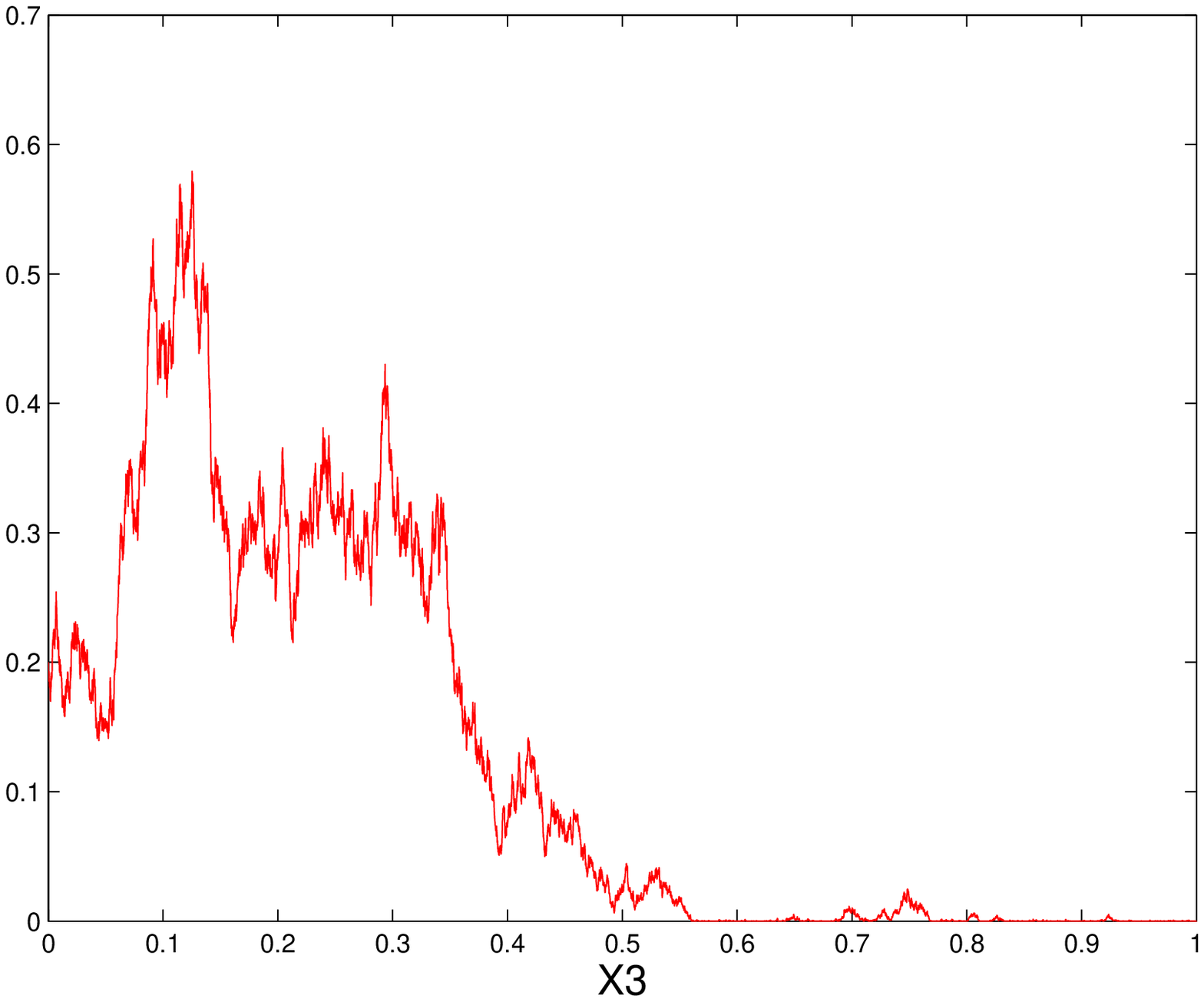}
\end{tabular}
\end{center}

The stochastic pendulum equation is considered in the following
manner. The dynamics of a non-dissipative classical pendulum of
the form $\ddot{\theta}(t)+2H\sin\theta(t)=0,$ can be expressed as
a system of stochastic differential equations expressed like
\begin{equation}\label{s5}
\left\{%
\begin{array}{ll}
    dx_1(t)=x_2(t)dt+\sqrt{x_1(t)}dW^1(t), &  \\
    dx_2(t)=-2H\sin(x_1(t))dt+\sqrt{x_2(t)}dW^2(t), &  \\
\end{array}%
\right.
\end{equation}and the Stratonovich stochastic pendulum equations
are
\begin{equation}\label{s5s}
\left\{%
\begin{array}{ll}
    dx_1(t)=(x_2(t)-\frac{1}{4})dt+\sqrt{x_1(t)}dW^1(t), &  \\
    dx_2(t)=-(2H\sin(x_1(t))+\frac{1}{4})dt+\sqrt{x_2(t)}dW^2(t), &  \\
\end{array}%
\right.
\end{equation}

For the probability density $p=u(x(t)),$ the corresponding
(It\^{o}) Fokker-Planck equation is given by
\begin{eqnarray*}
\frac{\partial}{\partial t}u(x(t))&=&-\frac{\partial}{\partial
x_1(t)}[x_2(t)u(x(t))]+\frac{\partial}{\partial x_2(t)}[2H
\sin(x_1(t))
u(x(t))]\\
&+&\frac{1}{2}\frac{\partial^2}{\partial
(x_1(t))^2}[x_1(t)u(x(t))]+\frac{1}{2}\frac{\partial^2}{\partial
(x_2(t))^2}[x_2(t)u(x(t))].
\end{eqnarray*}

Using stochastic Euler method on square root process, for initial
conditions $x_1(1)=1, \, x_2(1)=0.8$ we get the following graphics
for stochastic systems (\ref{s5}) and (\ref{s5s})
\begin{center}\begin{tabular}{ccc}
  \includegraphics[height=2in]{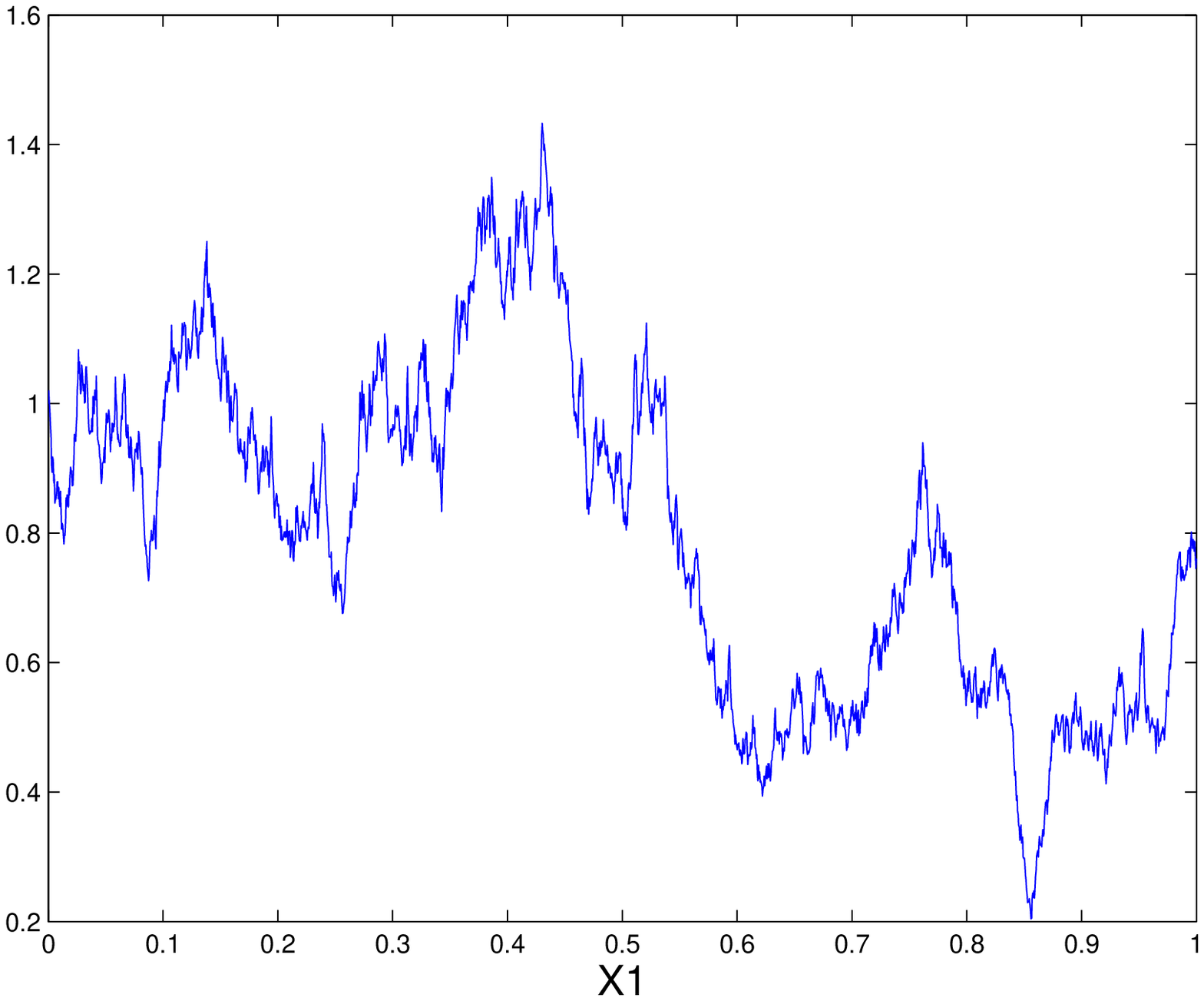}
  \includegraphics[height=2in]{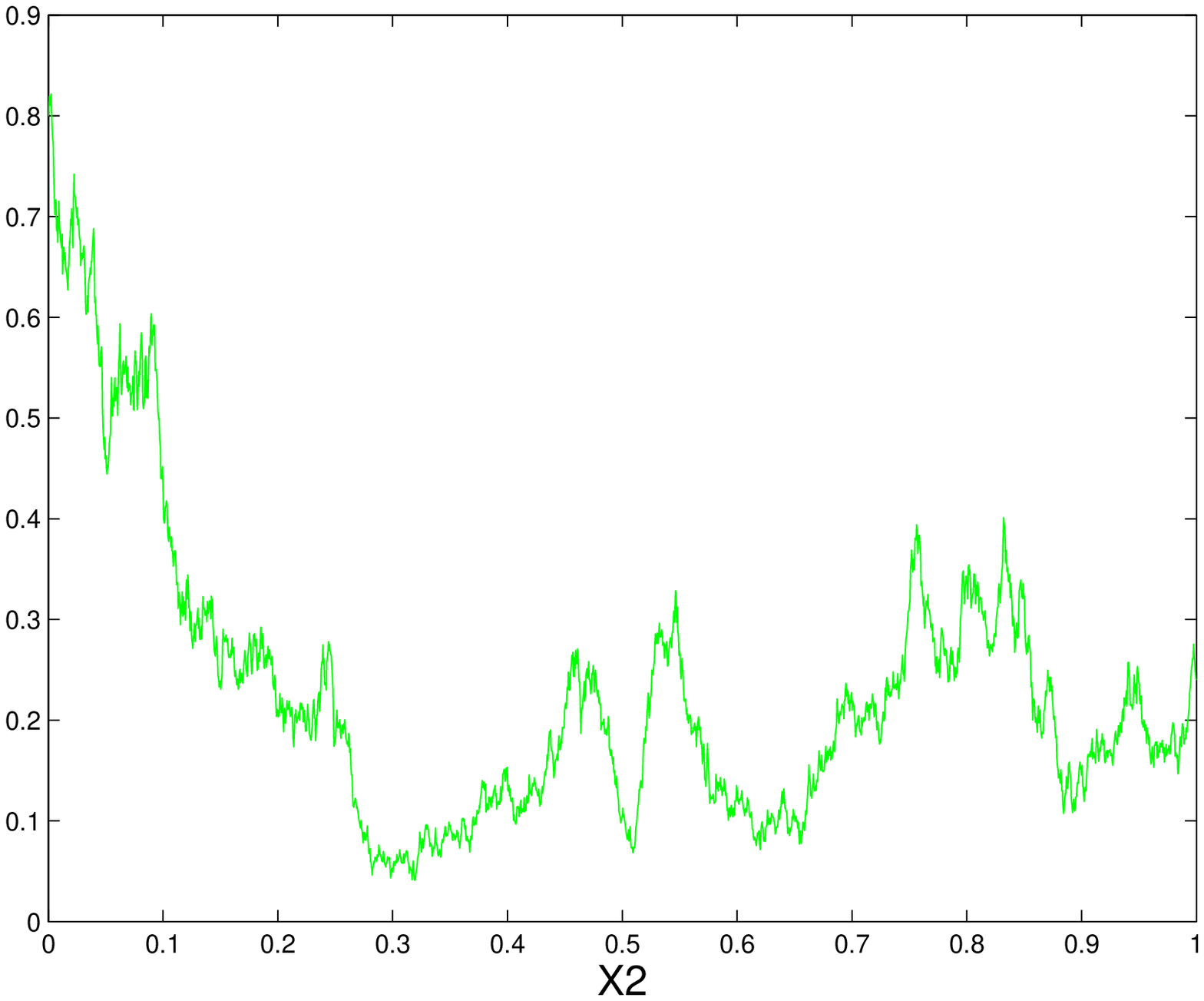}
\end{tabular}
\end{center}

\section{Conclusions}

In this paper we presented the Euler top system in $\mathbb{R}^3$
and the mathematical pendulum, but also the connections between
them: the existence of some applications that transform the
movement of a pendulum into a movement in $\mathbb{R}^3.$ That
means that the restriction of the Euler top system on a constant
level surface is the pendulum equation. This property is also true
in the case of Euler top system of differential equations with
delay argument, respectively mathematical pendulum with delay
argument, and in the case of fractional system of differential
equations, respectively fractional pendulum. We have also studied
the Euler top system and mathematical pendulum from the stochastic
point of view, using It\^o and Stratonovich integrals for a Wiener
process. Numerical simulations were done using Maple 12 and
Matlab. In the case of fractional Euler top system and fractional
pendulum we used the Adams-Moulton integration method for their
representation, and in the stochastic case we used the Milstein
scheme, that is a convergent numerical algorithm. In the future we
will study other aspects of these problems, such as stochastic
Lyapunov functions, stochastic Lyapunov exponents for determining
the stochastic stability in the equilibrium points  of a
considered system, classical, with delay of fractional.

\end{document}